\documentclass[a4paper, 11pt]{article}
\usepackage[utf8]{inputenc}
\usepackage{geometry}
\usepackage{graphicx}
\usepackage{xstring}
\usepackage{authblk} 
\usepackage{amsmath}
\usepackage{amssymb}
\usepackage{amsthm}
\usepackage{tcolorbox}
\usepackage{mathrsfs} 
\usepackage{stmaryrd} 

\usepackage{xcolor}
\usepackage{mathtools}

\usepackage{csquotes} 
\usepackage{hyperref}
\hypersetup{colorlinks,citecolor=blue,filecolor=blacklinkcolor=red, urlcolor=magenta,linkcolor=red}

\usepackage{tikz,tikz-cd} 

\newtheorem{defi}{Definition}
\newtheorem{lemm}{Lemma}
\newtheorem{coro}{Corollary}
\newtheorem{theo}{Theorem}
\newtheorem{prop}{Proposition}
\theoremstyle{remark}
\newtheorem{rema}{Remark}
\numberwithin{defi}{section}
\numberwithin{lemm}{section}
\numberwithin{coro}{section}
\numberwithin{theo}{section}
\numberwithin{prop}{section}
\numberwithin{rema}{section}
\numberwithin{equation}{section}

\usepackage[articlein=false,giveninits=true,style=ext-alphabetic,hyperref,backend=biber]{biblatex}
\AtEveryBibitem{%
  \clearfield{issn} 
  \clearfield{doi} 
  \ifentrytype{online}{}{
    \clearfield{url}
  }
}

\DeclareMathOperator{\sgn}{sgn}

\newcommand{\R}{\mathbb{R}}

\newcommand{\T}{\mathcal{T}}
\newcommand{\E}{\mathcal{E}}
\newcommand{\Q}{\mathcal{Q}}

\newcommand{\Ti}{\textnormal{Ti}}
\newcommand{\Spi}{\textnormal{Spi}}

\newcommand{\Weyl}[4]{W_{#1#2\phantom{#3}#4}^{\phantom{#1#2}#3}}

\newcommand*{\dscale}[1][2]{%
\IfEqCase{#1}{%
{-1}{|\bar{\lambda}|^{-\frac{1}{2}}}%
{1}{|\bar{\lambda}|^\frac{1}{2}}%
{2}{\bar{\lambda}}}}

\title{The Klein-Gordon equation in projective geometry} 
\author{Jack Borthwick\thanks{\href{mailto:jack.borthwick@u-paris.fr}{jack.borthwick@u-paris.fr}},\, Yannick Herfray\thanks{\href{mailto:yannick.herfray@univ-tours.fr}{yannick.herfray@univ-tours.fr}} }
\affil{*Université Paris Cité, Sorbonne Université, CNRS, IMJ-PRG, F-75013 Paris, France\\ $\dagger$Université de Tours, Université d’Orléans, CNRS, IDP UMR 7013, 37200 Tours, France\\
$\dagger$ Institut universitaire de France (IUF)}
\date{}

\begin{document}
\maketitle
{\let\thefootnote\relax\footnotetext{Keywords: projective geometry, Klein-Gordon equation, holography }}
\begin{abstract}
We discuss a geometric method by which Hörmander's result on the asymptotics of a Klein-Gordon field, on an asymptotically flat spacetime, can be related to some partial peeling, à la Penrose, along an extended projective boundary, embedded in a spacetime with one dimension higher.
\end{abstract}
\tableofcontents
\section{Introduction}

In 1987, Hörmander published a fascinating paper~\cite{hormanderRemarksKleinGordonEquation1987}\footnote{Slightly expanded upon in~\cite[Chapter VII]{Hormander:1997aa}} giving precise asymptotics for solutions to the Klein-Gordon equation:
\[(\Box -m^2)u =0,\]
with Cauchy data, $u=u_0$, $\partial_t u = u_1$ in the Schwartz class $\mathcal{S}(\R^{d+1})$.
A full Fourier transform of the equation leads to a natural Lorentz invariant splitting: $u=u_+ + u_-$ of the solution $u$ into positive and negative frequency parts, where the spectral projections are given by pseudo-differential operators:
\begin{equation}\label{eq:intro_projectors} u_+ = \frac{1}{2}\left( u+\frac{i}{m}\frac{\partial_t u}{\sqrt{-\Delta +1}}\right), \qquad u_-= \frac{1}{2}\left( u - \frac{i}{m}\frac{\partial_t u}{\sqrt{-\Delta+1}}\right). \end{equation}
If the solution is assumed to be real-valued then $u_+ = \bar{u}_-$; for the purpose of the present discussion, let us make this assumption and restrict our attention to $u_+$ which admits an integral formula:
\begin{equation}\label{eq:positive_frequency_integral_formula} u_+(t,x) =\frac{1}{(2\pi)^d} \int_{\R^d} \frac{a_+(\xi)}{\sqrt{1+|\xi|^2}} e^{-imt(\sqrt{1+|\xi|^2})+im\langle x,\xi \rangle}d\xi. \end{equation}
One can observe that: \[\xi \mapsto \frac{1}{m^{d-1}}\frac{a_+(\frac{\xi}{m})}{\sqrt{m^2+|\xi|^2}},\] is the Fourier transform of the initial value $u_+(0,\cdot)$. In his paper, Hörmander gives a very precise analysis of the asymptotic behaviour based on a formal application of the stationary phase method to this integral\footnote{Hörmander works with $m=1$, we have rescaled his formulae here for arbitrary mass $m$.} and shows that, up to a function in $\mathcal{S}(\R^{d+1})$, we have an asymptotic expansion: 
\begin{equation} \label{eq_asymptotic_expansion_KG} u_+ \sim (0^+ -i\varrho)^{\frac{d}{2}}e^{-\frac{im}{\varrho}}\sum_j \varrho^j w_j(t,x), \end{equation}
where $\varrho=\sgn(t)(t^2-|x|^2)^{-\frac{1}{2}}$ and,  importantly, $w_j(t,x)=w_j(1,\frac{x}{t})$.  The leading order term is given by:
\begin{equation}\label{eq:intro_leading_coeff} w_0(t,x)= \begin{cases} (2\pi)^{-\frac{d}{2}}a_+(x\varrho) & \textrm{if $t^2>|x|^2$}, \\ 0 & \textrm{if $t^2\leq |x|^2$}. \end{cases} \end{equation}
The dependence of the expansion coefficients $w_j$ on the ratio $\frac{x}{t}$ (rather than $x$, $t$ separately) is one of the first hints that projective geometry may play some role in a geometric interpretation of these asymptotics. 

\begin{figure}[h]\centering \includegraphics[scale=.7]{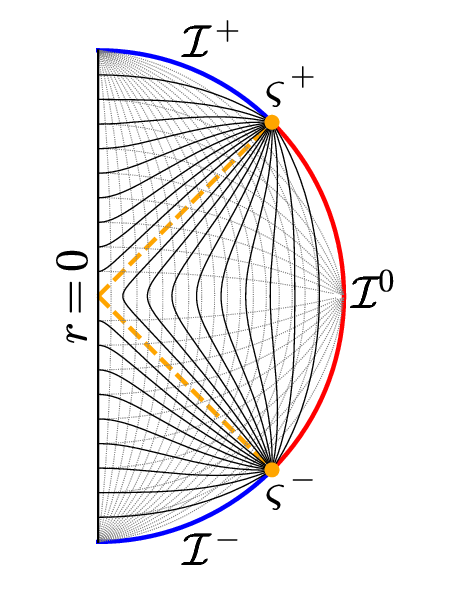}
\caption{\label{fig:CProjMinkowski}The projective compactification of Minkowski space time : $\mathcal{I^\pm}\simeq \mathbb{H}^d$ correspond to past and future time-like infinity, $\mathcal{I}^0\simeq\textrm{dS}^d$ is space-like infinity, and $\varsigma^\pm$, projective null infinity, are conformal spheres.  Reproduced from~\cite{Borthwick_Herfray_2025}.  }
\end{figure}

In the spirit of Penrose's conformal method~\cite{penrose_asymptotic_1963,penrose_zero_1965}, one might look to the \emph{projective compactification}~\cite{Cap:2014aa,Cap:2014ab,cap_projective_2016-1} as a means to tie this result to the asymptotic geometry of Minkowski spacetime. This compactification is instead based on the projective class of the Levi-Civita connection, made up of all torsion-free affine connections on $TM$ that share the same unparametrised geodesics, and constructs a boundary at infinity equipped with rich geometric structure, see Fig.~\ref{fig:CProjMinkowski}.

The function $\varrho=(t^2-|x|^2)^{\frac{1}{2}}$ is of specific relevance with respect to this construction: it constitutes a distinguished boundary defining function at (future/past) timelike infinity (denoted by $\mathcal{I}^\pm$) that puts the Minkowski metric in the (Beig-Schmidt~\cite{beig_einsteins_1982}) asymptotic form:
\begin{equation} \label{eq:MinkowskiBeigSchmidt} \varrho^2 \eta = -\frac{d\varrho^2}{\varrho^2} + h, \end{equation}
where $h$ is the hyperbolic metric. Introducing $y=\varrho x$, we can define a boundary coordinate chart $(\varrho, y)$ on the interior of the future/past lightcone, compatible with a hyperboloidal foliation such that $\rho=0$ is identified with timelike infinity. From here, it is clear that the leading order coefficient $w_0$ in Eq.~\eqref{eq:intro_leading_coeff}, can be interpreted as a function on $\mathcal{I}^\pm$.

In this geometric framework, it is extremely tempting to try to interpret the function $a_+$ as asymptotic data at timelike infinity, which enable to reconstruct the full solution via the integral formula~\eqref{eq:positive_frequency_integral_formula}, and to identify the \enquote{scattering operator}, that maps initial data on $\{t=0\}$ to the asymptotic data, with a warped Fourier transform (up to normalisation constants). However, the story is incomplete: as observed in~\cite{borthwickTiSpiCarrollian2025}, the limiting procedure producing the function $a_+$ is sensitive to the choice of boundary defining function $\varrho$ (preserving the asymptotic form~\eqref{eq:MinkowskiBeigSchmidt}). Although, it turns out, in a very controlled manner: only the lowest order change in the boundary defining function has an effect:
\[ \varrho \mapsto \varrho(1+ u(y)\varrho + O(\varrho^2)), \]
and obeys a very precise transformation rule:
\[ a_+ \mapsto a_+e^{imu}.\]
This can be interpreted geometrically in terms of an \emph{extended boundary} , which we refer to as $\Ti$, that is designed to encode geometrically this ambiguity, see Définition~\ref{DefinitionExtendedBoundary}. It is in fact a $\R$-principal bundle $\Ti \longrightarrow \mathcal{I}^\pm$: the transformation rule indicates that the \enquote{scattering data} is, rather than a function on its base, in fact, an equivariant function on $\Ti$ with respect to the one-dimensional representation\footnote{Observe here that the mass here is a choice of representation}: $(x,v) \mapsto e^{imx}v$ of $\R$. 
With this understanding of the scattering data, the interpretation of Eq.~\eqref{eq:positive_frequency_integral_formula} as a reconstruction formula from the boundary is even more striking. The integral is instead understood as being taken over a section of $\Ti$ (determined by the point of interest $(x,t)$) and Eq.~\eqref{eq:positive_frequency_integral_formula} is then assimilated with a Kirchoff-Adémar formula for massive fields, we refer the reader to~\cite[Proposition 4.1]{borthwickTiSpiCarrollian2025} for more details.

In the case of Minkowski spacetime, it is therefore apparent that $\Ti$ was, to some extent, a missing piece in the puzzle. Similar extended boundaries can be shown to exist for a class of manifolds defined to be asymptotically flat at time/spacelike infinity~\cite{Borthwick_Herfray_2025,borthwickTiSpiCarrollian2025} whose definition is inspired by the work of Ashtekar and Romano~\cite{ashtekar_spatial_1992}.   It is therefore natural to wonder what can be said about the Klein-Gordon equation in a curved setting. 

This paper is our first foray into this question with the goal to lay the foundations for a geometric framework for the study of the asymptotic behaviour of massive fields on Ricci-flat Lorentzian manifolds; emphasising the role of projective geometry throughout. For the sake of the present paper, which we believe should be thought of as a \enquote{proof of concept} rather than a definitive result, we return to the ambient space framework discussed in detail~\cite{RSTA20230042}, in which the Ricci-flat Lorentzian manifold $(\Q_{\mathcal{O}},h)$ is obtained as a quotient structure from a holonomy reduction of a projective geometry  on a higher-dimensional manifold $M \longrightarrow \Q$. These assumptions imply projective compactness of the metric $h$ in the sense of Čap-Gover~\cite{Cap:2014ab}, which is well-known to be too strong an assumption for many physically relevant cases, but which already showcases the important features of what can be recovered in an appropriately relaxed setting~\cite{capProjectiveInfinitiesBCalculus2026}.

The focus of our paper will be to exploit the ambient projective structure on $M$ to study the Klein-Gordon equation
\[ (h^{ab}\bar{\nabla}_a\bar{\nabla}_b - m^2)\varphi =0\,, \] on the open set $\Q_{\mathcal{O}}^+=\Q_{\mathcal{O}}\cup \{\bar{\sigma} >0\}$ of the base manifold, where we have denoted the Levi-Civita connection of $h$ by $\bar{\nabla}$.
As in the Penrose conformal method, this boundary behaviour will be understood to describe the asymptotic behaviour of the field.

After briefly recalling some elements of projective geometry in Section~\ref{sec:notations}, we complete, in Section~\ref{sec:ambient_geometry}, the results developed~\cite{RSTA20230042} by studying more precisely the geometry of the ambient space $M$ near \emph{its} boundary $\Sigma$ which, by construction, covers\footnote{see Fig.~\ref{fig:quotient_structure} below.} that of $\mathcal{Q}$. $\Sigma$ splits into curved orbits which are identified with the extended boundaries mentioned above. We recall that $M$ is equipped with a solution $\zeta^{ab}$ to the metrisability equation, the main result of this section is Proposition~\ref{prop:final_coordinate} in which we derive a normal form for  $\zeta^{ab}$ when expressed in distinguished projective scales linked to special local coordinate charts in a neighbourhood of the boundary. 

Section~\ref{sec:KG} is devoted to studying the Klein-Gordon equation on the base from the perspective of the ambient space. The first step is to find an appropriate projectively invariant differential operator that will, in some sense, englobe this equation. The natural candidate, the \enquote{projective Laplacian}, $H^{AB}D_AD_B$, does the trick, provided it is restricted to an appropriate subspace of projective densities that we call $\mathcal{E}_m(-\frac{n-1}{2})$.  Enforcing the condition $H^{AB}D_AD_B \phi =0$ on $\phi\in \mathcal{E}_m(-\frac{n-1}{2})$ shows that, when pulled-back to the base along any section of the (trivial) principal bundle $\mathcal{O}^+\longrightarrow \mathcal{Q}_{\mathcal{O}^+}$, we obtain a solution of the Klein-Gordon equation.

The structure of the space $\mathcal{E}_m(-\frac{n-1}{2})$, studied in Section~\ref{sec:m_density_space}, is particularly interesting: the \emph{ambient space fields} admit a natural decomposition $\phi=\phi_++\phi_-$ that recalls the positive/negative frequency decomposition~\eqref{eq:intro_projectors} of Hörmander. However, this decomposition is only recovered on the Klein-Gordon field on the base asymptotically. These considerations lead us to Theorems~\ref{thm_asymptotics} and~\ref{theo:asymptotics2}, that we would like to highlight: they constitute partial \emph{peeling}~\cite{penrose_zero_1965,mason_nicolas_2009} results, within our framework, that relate the leading order asymptotic behaviour in Eq.~\eqref{eq_asymptotic_expansion_KG} to regularity extensions of the ambient fields $\phi_\pm$ up to the projective boundary. These results are only partial as we do not consider, in this paper, the question of reconstructing the ambient field from the base field. Indeed, our positive and negative frequency decomposition is only recovered asymptotically (see Lemma~\eqref{lemm:recover_decomp}), and hence we believe that reconstructing the solution from this data is a question of a different analytical nature and choose here to emphasise the geometric and structural aspects, postponing these considerations to future work.

In Section~\ref{sec:kg_boundary_calculus}, we show that the geometric framework also provides a \emph{boundary calculus} for the Klein-Gordon equation. This result is to be related to the works of~\cite{Gover:2014aa,gover_poincare-einstein_2015} on Poincaré-Einstein manifolds in conformal geometry (see also~\cite{Curry_Gover_2018} for a pedagogical introduction). They exploit a natural $\mathfrak{sl}_2$-algebra associated with the Degenerate Laplacian that can be used to solve boundary value problems, \emph{formally} off the boundary by a natural iterative procedure. The algebra here is \emph{not} $\mathfrak{sl}_2$ nor its Inönü contraction: it is instead the Heisenberg algebra.

Finally, we conclude the article in Section~\ref{sec:model} by studying the important model case.

  \medskip

We would like to mention that, in the course of writing this work, we became aware of independent work~\cite{markov2026} who have obtained similar results by completely different methods, we look forward to discussing with these authors.
\section*{Acknowledgements}
Although already two years ago, the first ideas for this work were discussed on the tram in Vienna during the workshop \enquote{Carrollian physics and holography}(2024); the authors are thankful once more to the Erwan Schrödinger for hosting them during this period. This work is supported by the ANR project~\emph{ANR-25-CE40-4883 ScattHoloGR}. JB would like to thank Andréas Čap for several interesting discussions on the results in this work.

\section*{AI disclosure}
All the ideas of the paper are those of the authors and no genAI technology was used at any point in the writing process of this paper. 
Rights to engage with this manuscript using generative AI technology are \textbf{not} granted. JB refuses permission to use this manuscript in training AI models. 
\section{Notations}\label{sec:notations}
It should be understood that all affine connections on the tangent bundle are torsion-free.
We make extensive use of the Penrose abstract index notation~\cite{penrose_spinors_1984}. In index notation, tractor indices (see below) will be denoted by capital Latin letters, $A,B,C\dots,$. 
\subsection{Projective geometry and tractor calculus}
Let $M$ be a manifold of dimension $n+1$ endowed with an equivalence class of projectively equivalent torsion-free connections $[\nabla]$. It is a standard result\footnote{See, for instance,~\cite[Proposition 7.2]{Kobayashi:1995aa}} that two torsion-free affine connections $\hat{\nabla}$ and $\nabla$ are projectively equivalent if and only if there is a 1-form $\Upsilon_a$ such that for any vector field $\xi^b$:
\[\hat{\nabla}_a\xi^b=\nabla_a\xi^b + \Upsilon_a\xi^b +\Upsilon_c\xi^c\delta_a^b. \]
We will use $\hat{\nabla}=\nabla + \Upsilon$ as shorthand for this condition.

We denote by $\mathcal{E}(w)$ the associated vector bundle to the frame bundle $P^1(M)$ determined by the $1$-dimensional representation of $GL_{n+1}(\R)$:
\[ A \mapsto \lvert \det {A}\rvert ^{\frac{w}{n+2}}.\]
Sections of this bundle will be called \emph{projective densities} of weight $w$. Their definition is so that for any $\sigma \in \Gamma(\mathcal{E}(w))$ under the projective change of connection $\hat{\nabla} = \nabla + \Upsilon$, \[\hat{\nabla}\sigma = \nabla \sigma + w\Upsilon \sigma. \]
For any vector bundle $\mathcal{B}$ with base $M$, we will write: $\mathcal{B}(w)=\mathcal{B}\otimes \mathcal{E}(w)$ and we will say that sections of $\mathcal{B}(w)$ are of weight $w$.

Let $G=\textnormal{PGL}(n+2,\R)$ and $H=(\R^{n+1})^*\rtimes GL(n+1,\R)$ the isotropy subgroup of a given line. The class $[\nabla]$ determines a reduction of the second-order frame bundle $P^2(M)$ to a $H$-principle bundle $P$, on which there is a unique \emph{normal} Cartan connection $\omega$~\cite{Cartan:1924aa, Kobayashi:1995aa}; this is a normal projective geometry modeled on projective geometry as described in~\cite{Sharpe:1997aa}.

The \emph{standard (projective) tractor bundle}, $\T$ can be defined as the associated vector bundle to $P$ determined by the restriction to $H$ of the $(n+2)$-dimensional representation of $G$ given by:  $A\mapsto \lvert\det A\rvert^{-\frac{1}{n+2}}A$. The Cartan connection $\omega$ induces a linear connection $\nabla^{\mathcal{T}}$ on $\mathcal{T}$.

The bundle $\mathcal{T}$ fits into a canonical short exact sequence of vector bundles:
\begin{center}
\begin{tikzcd} 0 \arrow[r] & \mathcal{E}(-1) \arrow[r,"X"] & \arrow[l,dotted,bend left = 45, "Y"] \mathcal{T} \arrow[r,"Z"] & TM(-1) \arrow[r] \arrow[l,dotted,bend left=45,"W"] &0. \end{tikzcd}
\end{center}
A choice of connection in the projective class $[\nabla]$ provides a non-canonical isomorphism $\mathcal{T}\simeq \mathcal{E}(-1)\oplus TM(-1)$, described above by the non-canonical maps $Y$ and $W$. When thinking of these maps as sections of $\mathcal{T}^*\otimes\mathcal{E}(-1)$ and $T^*M(1)\otimes\mathcal{T}$ respectively, they transform under change of connection $\hat{\nabla}=\nabla + \Upsilon$ according to:
\begin{equation}\hat{Y}_A=Y_A - \Upsilon_aZ^a_A, \quad \widehat{W}^A_a=W^A_a + \Upsilon_a X^A. \end{equation}
Fixing a choice of connection $\nabla$ in the projective class, sections of the standard tractor bundle will be written:
\[ T^A=\xi^aW_a^A + \rho X^A,\]
similar expressions will be used with sections of tensor powers of the tractor bundle.

The dual tractor bundle $\T^*$ is usually identified with the first jet prolongation $J^1\E(1)$ of the weight $1$ projective density bundle~\cite{Bailey:1994aa}. This description has the advantage of being more direct than the principal bundle approach outlined above, and is independent of the projective structure. In addition, it allows for straightforward comparisons of tractors on different manifolds; in particular, tractors in the \enquote{bulk} and on distinguished hypersurfaces.

We recall from~\cite{Bailey:1994aa} the standard decomposition of the curvature tensor $R_{ab\phantom{c}d}^{\phantom{ab}c}$:
\[ R_{ab\phantom{c}d}^{\phantom{ab}c}=W_{ab\phantom{c}d}^{\phantom{ab}c} + 2\delta^c_{[a}P_{b]d} + \beta_{ab}\delta^c_d,\]
where the Weyl tensor $W$ is trace-free and $P_{ab}$ is the projective Schouten tensor. It is easily verified that:
\[nP_{ab}= R_{ab}+\beta_{ab},\quad \beta_{ab}=-\frac{2}{n+2}R_{[ab]}=-2P_{[ab]},\]
where $R_{ab}=R_{da\phantom{d}b}^{\phantom{da}d}$ is the Ricci tensor. The tensor $\beta_{ab}$ can be thought of as the curvature tensor of the density bundle $\mathcal{E}(1)$.

A connection $\nabla$ is said to be \emph{special} if it preserves a nowhere vanishing density $\sigma$. In this case $\beta_{ab}=0$ (in particular all density bundles are flat) and we say that $\nabla$ is the \emph{scale} determined by $\sigma$. Even when it does not preserve a nowhere vanishing density, we will refer to a choice of connection in the class as a choice of scale.

Having fixed a projective connection $\nabla$ in the class, and split the short-exact sequence, the action of the normal Cartan connection can be summarised conveniently by the relations:
\begin{equation}\label{TractorConnection} \nabla_a X^A= W^A_a, \quad \nabla_a W^B_b=-P_{ab}X^B,\quad \nabla_a Y_A=P_{ab}Z^b_A, \quad \nabla_a Z^b_B=-\delta_{a}^bY_B. \end{equation}

The tractor curvature, defined by, the identity $\Omega_{ab\phantom{C}D}^{\phantom{ab}C}T^D= 2 \nabla_{[a}\nabla_{b]}T^C$, is known to be given in an arbitrary scale by:
\begin{equation}\label{TractorCurvature} \Omega_{ab\phantom{C}D}^{\phantom{ab}C}=\Weyl{a}{b}{c}{d}W^C_cZ^d_D -Y_{abd}X_CZ^d_D \end{equation}
Where: $Y_{abd}=2\nabla_{[a}P_{b]d}$ is the projective Cotton tensor.

\subsection{Metric projective structures}
A projective class $[\nabla]$ is said to be \emph{metric}, if it admits a solution $\zeta^{ab} \in \Gamma(S^2(T^*M)(-2))$ to the Mike\v{s}-Sinjukov metrisability equation~\cite{Mikes1996,Sinjukov:1979aa}:
\[ \textnormal{trace-free}\left( \nabla_c\zeta^{ab} \right)=0, \]
or equivalently $\nabla_c \zeta^{ab} = -2 \delta^{(a}_c \lambda^{b)}$ for some $\lambda^a$.
Solutions of this equation have been shown~\cite{Eastwood:2008aa,Flood:2018aa} to be in one-to-one correspondence with tractors $H^{AB}$ that satisfy:
\begin{equation} \label{TractorMetrisability}\nabla_c H^{AB} + \frac{2}{n+1}X^{(A}W_{cE\phantom{B}F}^{\phantom{cE}B)}H^{EF}=0, \end{equation}
where $W_{cE\phantom{B}F}^{\phantom{cE}B}=Z^c_C\Omega_{ce\phantom{B}F}^{\phantom{ce}B}$ and $\Omega_{ce\phantom{B}F}^{\phantom{ce}B}$ is the tractor curvature.
In this case $H^{AB}$ is of the form:
\[H^{AB}= \zeta^{ab}W^A_aW^B_b + 2\lambda^aX^{(A}W^{B)}_b+ \tau X^AX^B.\]
with 
\begin{align}\label{H decomposition}
\lambda^a &= - \frac{1}{n+2}\nabla_c\zeta^{ca},& \tau&=  \frac{\nabla_a\nabla_b \zeta^{ab}}{(n+1)(n+2)}+\frac{1}{n+1}P_{ab}\zeta^{ab}.
\end{align}

If $\nabla_c H^{AB}=0$, then it can be shown that Eq.~\eqref{TractorMetrisability} is automatically satisfied; these solutions are said to be \emph{normal}. This is equivalent to
\begin{align}\label{DH=0}
\nabla_c \zeta^{ab} + 2 \delta^{(a}_c \lambda^{b)} &=0, &
\nabla_c \lambda^a + \delta^a_c \tau - P_{cb}\zeta^{ba} &=0,
\end{align}
(and imply \eqref{H decomposition}). Normal solutions are known to correspond to Einstein metrics~\cite{Cap:2014aa,Armstrong2008}. 
\subsection{The Thomas $D$-operator}
For any weighted tractor bundle, $T^\mathscr{A}(w)$, where $\mathscr{A}$ is an arbitrary set of abstract indices, the normal Cartan connection induces a natural projectively invariant differential operator:
\[ D: T^\mathscr{A}(w) \rightarrow \T_A\otimes T^\mathscr{A}(w-1)\]
known as the Thomas $D$-operator and defined by:
\begin{equation}\label{Definition: Thomas operator}
D_A t^\mathscr{A} = w\,t^\mathscr{A}\,Y_A + \nabla_a t^\mathscr{A}\,Z^a_A.
\end{equation} 
$D$ behaves much like a covariant derivative with a cotractor index, in particular, \emph{it satisfies the Leibniz rule} and is iterable; this can be useful in computations. Also observe that one has
\begin{equation}\label{DderivativeX} D_AX^B = \delta_{\,\,A}^{B}. \end{equation}

\section{Ambient geometry near the projective boundary}\label{sec:ambient_geometry}
\subsection{Setup}
Since the present work builds heavily on~\cite{RSTA20230042} we will follow, as closely as possible, the notation used there for ease of reference; we began by recalling the initial setup.
Throughout $\widetilde{M}$ will denote a compact manifold $\widetilde{M}$ of dimension $n+1$, where: $n=d+1$. We assume $\widetilde{M}$ is equipped with a projective structure $[\nabla]$ of which the corresponding normal Cartan connection admits a holonomy reduction given by the existence of two parallel tractors, $H^{AB}, I^A$, i.e. such that:
\begin{equation}\label{eq:IandHparallel} \nabla_c I^A= 0, \qquad \nabla_c H^{AB}=0.\end{equation}We assume that $H^{AB}$ is everywhere non-degenerate with signature $(d,2)$ and inverse $\Phi_{AB}$ \emph{and}
\[ \Phi_{AB}I^AI^B=0.\]
The construction comes with two distinguished projective densities, of weight $1$ and $2$ respectively given by:
\begin{equation}\label{eq:Inull} \sigma=\Phi_{AB}I^AX^B, \qquad \lambda=\Phi_{AB}X^AX^B,\end{equation}
a distinguished weighted vector field (weight $-1$):
\[ I^AZ_A^a= n^a,\]
and a weight $-2$ bilinear form: \[ \zeta^{ab}=H^{AB}Z_A^aZ_B^b.\]

It can be useful to have in mind the following identities resulting form the differential assumptions~\eqref{eq:IandHparallel} :
\begin{equation}\label{eq:CsqParallel} D_A\sigma = \Phi_{AB}I^B, \qquad \frac{1}{2}D_A\lambda = \Phi_{AB}X^A,\qquad  \Phi_{AB}=\frac{1}{2}D_AD_B\lambda, \end{equation}
where $D$ is the Thomas $D$-operator.

We recall from~\cite[Theorem 3.1]{RSTA20230042} the description of the curved orbit decomposition of $\widetilde{M}$:
\[ \widetilde{M}= \overbrace{\mathcal{O}\sqcup \underbrace{{\Sigma_+}\sqcup{\Sigma_0}\sqcup{\Sigma_-}}_{\Sigma}}^M \sqcup\mathcal{Z}(n), \]
where the open orbit $\mathcal{O}\subset M$ satisfies
\begin{equation}
	\mathcal{O} = \{\sigma \neq 0\}
\end{equation}
and its boundary, $\Sigma \subset M$, is the union of  \[\Sigma_{\pm}=\mathcal{Z}(\sigma)\cap\{\pm \lambda >0\}, \quad \Sigma_0=\mathcal{Z}(\sigma)\cap\{\lambda=0\}\,.\]

It is important to emphasise the different roles played by the densities $\sigma$ and $\lambda$: $\sigma$ describes the location of the boundary $\Sigma$, whereas the role of $\lambda$ is to describe its finer decomposition. For our purposes, our main focus will be the components\footnote{One generically expects some level of singularity  along $\Sigma_0$, induced by gravitational radiations. We will not try to address these here.} $\Sigma^\pm$ which are the bundles $\Ti$ and $\Spi$, the first of which is expected~\cite{borthwickTiSpiCarrollian2025} to carry the asymptotic data for massive fields.

As observed in~\cite[Theorem 3.2]{RSTA20230042} the flow of $n^a$ preserves the above decomposition and we assume the leaf space $\Q$ of $M$ under this flow is a smooth manifold; this is always the case on the open orbit $\mathcal{O}$.

We denote by $\mathcal{Q}_{\mathcal{O}}$ the image of $\mathcal{O}$ under the canonical projection $\pi : M \rightarrow \Q$, it is equipped with a projectively compact Ricci flat Lorentzian\footnote{signature $(d,1)$} metric $h$ (induced from $g^{ab}=\sigma^{2}\zeta^{ab}$ on $\mathcal{O}$) see~\cite[Proposition 4.7]{RSTA20230042} , whose projective boundary is given by $\pi(\Sigma)$ (see~Fig.~\ref{fig:quotient_structure}).
\begin{figure}[h!]
\centering
\begin{tikzcd}[column sep=2mm]
M\arrow[d] &=&\mathcal{O}\arrow[d] &\sqcup& \overbrace{ \left(\Sigma^- \sqcup \Sigma^0 \sqcup \Sigma^+\right)}^\Sigma\arrow[d]\\
\Q  &\simeq& \Q_\mathcal{O} &\sqcup&  \underbrace{\left(  \mathcal{H}^- \sqcup  \mathcal{H}^0 \sqcup   \mathcal{H}^+\right)}_{\mathcal{H}}
\end{tikzcd}
\caption{\label{fig:quotient_structure} Links between the ambient and quotient structures}
\end{figure}

Additionally, not only $\mathcal{O}$ is a principal $\R$-bundle over $\Q_{\mathcal{O}}$ but it comes with canonical global trivialisation
\begin{equation}
	\mathcal{O} = \mathbb{R} \times \Q_{\mathcal{O}}
\end{equation}
 given by the map: $x \mapsto (\frac{1}{2}l(x),\pi(x))$ where we define:
\begin{equation} l = \sigma^{-2}\lambda.\end{equation} 
Similarly, the curved orbits of the boundary $\Sigma^\pm$ (resp. $\Sigma^0$), which are strong Carrollian geometries (resp. strong conformal Carrollian geometries)\cite{duval_carroll_2014,duval_conformal_2014,Herfray:2021qmp} (see e.g. \cite{Bekaert:2026cvx} for a pedagogical introduction), are principal $\R$-bundle over the curved orbits $\mathcal{H}^{\pm}$ (resp $\mathcal{H}^{0}$ ) of the projective compactification of $\Q_{\mathcal{O}}$.

Importantly, these structures can be recovered from the data of the projective compactification on the base manifold $\mathcal{Q}$ alone~\cite{Borthwick_Herfray_2025}. We recall the intrinsic definitions of these extended boundaries:
\begin{defi}\label{DefinitionExtendedBoundary}
The \textbf{extended boundaries} $\Spi/ \Ti$ are defined to be the pullback to $\mathcal{H}^{\pm}$ of $i^*J^1\mathcal{E}(1)_{\mathcal{Q}}=i^*\T^*_{\Q}$ (viewed as a bundle with base $i^*J^1\mathcal{E}(1)_{\Q}/J^1\bar{\sigma} \simeq \T^*_{\mathcal{H}}$) along the Thomas derivative $D_A \dscale[1]$: 
\begin{center}
\begin{tikzcd}
\Sigma^\pm \quad  \arrow[d, "\pi"] \arrow[r,dashed] &\quad  i^*\T^*_\Q \arrow[d,"\textnormal{Can. Proj.}"] \\[0.7em] \mathcal{H}_{\pm}\quad  \arrow[r, "D_A{\dscale[1]}"]&\quad \T^*_\mathcal{H}\simeq \T^*_\Q/I
\end{tikzcd}
\end{center}
\end{defi}
Here $\bar{\sigma}$ is the boundary defining density of the projective compactification and $\bar{\lambda}$ is the density \emph{along the projective boundary $\mathcal{H}$} describing its decomposition. They are related to the ambient densities $\sigma$  and $\lambda$ by: 
\[ \sigma = \pi^*\bar{\sigma},\qquad \lambda = \pi^*\bar{\lambda},\]
where we emphasise that the second inequality only holds on $\Sigma$ (see Eq.~\eqref{eq:projective_lie_derivative_lambda} below).

The sections of these extended boundaries have an interpretation in terms of an important class of scales on the base:
\begin{theo}[Theorem 2.1~\cite{Borthwick_Herfray_2025} reformulated from~\cite{Cap:2014ab}]\label{thm:cap_gover_scale} Let $\Q=\Q_{\mathcal{O}}^+ \cup \mathcal{H}^\pm \cup \mathcal{H}^0$ be a smooth projective compactification of a Ricci flat metric $h$ on $\Q_{\mathcal{O}}^+$, $\bar{\zeta}^{ab}$ the corresponding solution to the metrisability equation and $\bar{\sigma}$ the canonical boundary defining density.
A \underline{special} scale $\nabla \in [\nabla^h]$, defined on a neighbourhood $U$ of a point $x_0 \in \mathcal{H}^\pm$, is such that 
\[ \bar{N}^a: =\bar{\sigma}^{-2}\bar{\zeta}^{ab}\nabla_b \bar{\sigma}, \]
extends smoothly to $\mathcal{H}^\pm\cap U$ if and only if it preserves a density $\tau$ satisfying
\[ i^*\tau = |\bar{\lambda}|^\frac{1}{2}.\] 

When evaluated in such scales the metric can then be rewritten as
\begin{equation}\label{Projective compactification review: asymptotic form of the metric} \bar{\zeta}_{ab}= \frac{\nu^{-1}}{\bar{\sigma}^2}\nabla_a\bar{\sigma} \nabla_b\bar{\sigma} + \bar{q}_{ab},
		\end{equation}
		for some finite tensorial density $\bar{q}_{ab}$ and $\nu=\bar{N}^a\nabla_a\bar{\sigma}$.
\end{theo}

\subsection{Boundary scales}\label{sec:boundary_scales}
We remind the reader that for any $w\in \mathbb{R}$, there is a canonical identification (\cite[Proposition 5.1]{RSTA20230042}): \[\pi^*\mathcal{E}_{\Q}(w)\simeq \mathcal{E}_M(w).\] An important part of the dictionary between structures and fields on $M$ and those on $\Q$ is that pullback densities, which we will refer to as \emph{adapted scales}, (see also~\cite[Definition 5.5]{RSTA20230042}) lie exactly in the kernel of the projectively invariant operator: \[ D_n: \mathcal{E}(w) \rightarrow \mathcal{E}(w-1)\,,\] defined by: \[ s \mapsto I^AD_A s.\]
It follows from Eqs.~\eqref{eq:Inull}, \eqref{eq:CsqParallel} that $\sigma$ is such a density which justifies the equality $\sigma = \pi^*\bar{\sigma}$. However, $\lambda$ is \emph{not} globally an adapted scale because: 
\begin{equation}\label{eq:projective_lie_derivative_lambda}  D_n\lambda =2\sigma.\end{equation} 
Whilst for the quotient structure this scale is geometrically relevant only along $\{\sigma=0\}$, it is nevertheless of greater significance with respect to the ambient projective structure on $M$. Indeed, the connection $\hat{\nabla}$ in the projective class defined by the property $\hat{\nabla}_a\lambda=0$ is the Levi-Civita connection of the metric: \[ \hat{g}_{ab}=\sgn(\lambda)\lambda^{-1}\zeta_{ab},\]
on the ambient space. Additionally, $\lambda$ does not vanish on $\Sigma^\pm$ and therefore can be used to study the local geometry in the ambient near these regions.
Using our assumptions on $I^A$ and $H^{AB}$, we have the following property in the scale $\hat{\nabla}$ determined by $\lambda$.
\begin{lemm}\label{lemm:properties_scale_lambda}
\[ \begin{gathered}\hat{\nabla}_b n^a=-\sigma\lambda^{-1}\delta_b^a,\qquad \zeta^{ab}\hat{\nabla}_b\sigma = n^a, \qquad I^A= n^a\hat{W}_a^A + \sigma\lambda^{-1}X^A \vspace{1mm}\\ 
H^{AB}=\zeta^{ab}\hat{W}_a^A\hat{W}_b^B + \lambda^{-1} X^AX^B.
 \end{gathered} \]
\end{lemm}
\begin{proof}The expression for $H^{AB}$ follows immediately from the fact that $\hat{\nabla}_c\zeta^{ab}=0$ and that its inverse $\Phi_{AB}=\frac{1}{2}D_AD_B\lambda.$ The expression for $I^A$ is a direct application of~\cite[Lemma 3.8 (3.15)]{RSTA20230042} to $\hat{\nabla}$.
Now since $D_A\sigma = \Phi_{AB}I^B$ we deduce that: $H^{AB}D_B\sigma = I^A$, from which the second equation follows.
The first equation is just~\cite[(3.2)]{RSTA20230042} combined with the expression for $I^A$.
\end{proof}

The fact that the scale $\lambda$ itself is not well-aligned with the flow lines of $n^a$, means that it will not provide us with a clean dictionary with the structure on the leaf space $\Q$.  Our present goal is to exploit the properties of $\lambda$ to find a class of adapted scales that lead to a normal form for the ambient weighted metric $\zeta^{ab}$. As one might expect, they will turn out to be exactly the pullback of scales satisfying the hypotheses of Theorem~\ref{thm:cap_gover_scale}.

We begin by introducing an arbitrary local coordinate system: \begin{equation}
	(u,x) \in \mathbb{R} \times \mathcal{Q}
\end{equation}
on $M\cap \{\lambda\neq 0\}$ such that: \[\hat{N}^a=\sqrt{|\lambda|}n^a=\partial^a_u.\]
Just as the ratio between the scales $l=\lambda\sigma^{-2}$ was important in studying the open orbit $\mathcal{O}=M\cap\{\sigma\neq 0\}$, we consider here the ratio:
\begin{equation}
	s= \frac{1}{\sqrt{l}} = \sigma|\lambda|^{-\frac{1}{2}},
\end{equation}
by definition it is a boundary defining function for $\Sigma^\pm \subset M$. %

Trivialising the equations~\eqref{lemm:properties_scale_lambda} in the scale $\lambda$ we have:
\[\hat{g}^{ab}\hat{\nabla}_b s = \hat{N}^a,\qquad \hat{\nabla}_b \hat{N}^a=-s\sgn(\lambda)\delta^{a}_b\]
Since \[n^an_a=-\sigma^2\lambda^{-1} \Rightarrow \hat{N}^a\hat{N}_{a}=-\sgn(\lambda)s^2,\]
it follows that:
\[ \hat{N}^a\hat{\nabla}_as=-\sgn(\lambda)s^2.\]
Solving this differential equation leads to a relationship of the following type between $s$ and $u$:
\[ s(u,x)=\frac{\varrho}{\sgn(\lambda)\varrho(u-u_0) +1}, \]
for some function $\varrho$ such that $\mathcal{L}_n\varrho =0$. 
\begin{rema}
This result shows that the principle bundle structures on the interior and on the boundary are incompatible.
\end{rema}
Shifting $u$, we can always assume $u_0=0$.

 We shall work in a region near $\Sigma^{\pm} = \{\varrho=0\}\cap \{\pm \lambda >0\}$, and define the following quantity:
 \[ \tau= \frac{1}{(1+\sgn(\lambda)\varrho u)^2}\lambda=\frac{s^2}{\varrho^2}\lambda.\]
 Note, in particular, that along $\Sigma^{\pm} = \{\rho=0\}\cap \{\pm\lambda >0\}$, one has  $\tau = \lambda$.
 
 \begin{lemm}
 $\tau\in \Gamma\mathcal{E}(2)$ is an adapted scale:
 \[ D_n\tau =0.\]
 In particular, if $\bar{\nabla}$ is the scale preserved by $\tau$ then:
\[ \bar{\nabla}_a n^b =0, \qquad n^a\bar{\nabla}_a\sigma =0. \]
 \end{lemm}
\begin{proof} \[\begin{aligned} D_n\tau = \hat{\nabla}_n (\frac{s^2}{\varrho^2})\lambda +2\sigma\lambda^{-1}\frac{s^2}{\varrho^2}\lambda&=2\frac{s}{\varrho}\left(\frac{n^a\nabla_a s}{\varrho} \right)\lambda+2\sigma\frac{s^2}{\varrho^2},\\&=2\frac{s}{\varrho^2}(-s\sigma\lambda^{-1})\lambda+2\sigma\frac{s^2}{\varrho^2},\\&=0. \end{aligned}\]
We record from this computation that:
\begin{equation} n^a\hat{\nabla}_a \tau= -2\frac{s^3}{\varrho^2}|\lambda|^{\frac{1}{2}}=-2\sigma \tau \lambda^{-1}, \end{equation}
\end{proof}

Let us now consider the properties of the scale $\bar{\nabla}$ that preserves $\tau$, we define the corresponding tensor fields: \begin{equation}\label{eq:tensors_scale_tau}\bar{g}_{ab}=|\tau|^{-1}\zeta_{ab},\qquad \nu = \sigma^{-2}\zeta^{ab}\bar{\nabla}_a\sigma\bar{\nabla}_b\sigma,\qquad t^b=\sigma^{-2}\zeta^{ab}\bar{\nabla}_a\sigma.  \end{equation}
We begin with:
\begin{lemm}\label{lemm:properties_nu_t}
\[\zeta_{ab}t^an^b=0, \qquad t^a\nabla_a u = \sigma^{-2}\sgn(\lambda)|\tau|^{\frac{1}{2}}\nu, \qquad D_n\nu =0.\]
\end{lemm}
\begin{proof}
The first equality is due to the fact that $n^a\bar{\nabla}_a\sigma =0$, the second then follows from:
\[\nabla_a s= \frac{s^2}{\varrho^2}\nabla_a\varrho - \sgn(\lambda)s^2\nabla_a u, \]
which implies that:
\[\zeta_{ab}n^b=\hat{\nabla}_a \sigma=|\tau|^{\frac{1}{2}}|\lambda|^{-\frac{1}{2}}\bar{\nabla}_a\sigma-\sigma^2\sgn(\lambda)|\lambda|^{-\frac{1}{2}}\nabla_a u. \]
From which we obtain the second equality by contracting with $t^a$.
For the third, we continue to work in the scale $\bar{\nabla}$; since it is special: $D_n \nu =\bar{\nabla}_n \nu$.
Using that $n^a\bar{\nabla}_a\sigma =0$ and $\bar{\nabla}_a\bar{\nabla}_b\sigma =-\bar{P}_{ab}\sigma$, we obtain that:
\[\begin{aligned} \bar{\nabla}_c\nu &=2\sigma^{-1}\nu \bar{\nabla}_c \sigma - 2 \sigma^{-1}\nu\bar{\nabla}_c\sigma-2\zeta^{ab}\sigma^{-1}\bar{P}_{ac}\bar{\nabla}_b\sigma&=-2\sigma t^a\bar{P}_{ac}. \end{aligned} \]
The last claim is then a consequence of the fact that (see~\cite[Lemma 5.6]{RSTA20230042}):
\[n^c\bar{P}_{ca}=0\] because $\tau$ is an adapted scale.
\end{proof}
Since, by construction, $\tau \to \lambda$ when $\varrho \to 0$, and $D_n\tau=0$ we have the following:
\begin{lemm}\label{lemm:properties_scale}
The scale induced by $\tau$ on the quotient $\mathcal{Q}$ satisfies the hypotheses of Theorem~\ref{thm:cap_gover_scale}. In particular:\begin{enumerate} \item $\nu$ extends to the boundary (with non-vanishing boundary value), \item 
\[ t^a= \sigma^{-2}\sgn(\lambda)|\tau|^\frac{1}{2}|\lambda|^{\frac{1}{2}}\nu n^a + \bar{t}^a, \]
with $\bar{t}^a\nabla_a u =0$ and $\bar{t}^a$ has a finite non-vanishing extension to the boundary $\sigma=0$.
\end{enumerate}
\end{lemm}

We now consider the backwards construction: instead of choosing a coordinate system, one can also work from a choice of adapted scale:
\begin{prop}\label{prop:coordinate_step1}
Let $\tau$ be an arbitrary adapted scale ($D_n\tau =0$), such that $\iota^*\tau=\lambda$, where $\iota$ denotes the canonical inclusion of $\{\sigma=0\}$ in $M$, then one can always choose a coordinate system $(u,x)$ such that $\partial^a_u = \hat{N}^a$ and, if $\varrho=\sigma|\tau|^{-\frac{1}{2}}$, then:
\[ s=\sigma|\lambda|^{-\frac{1}{2}} = \frac{\varrho}{1+\sgn(\lambda)\varrho u}.  \]
\end{prop}
\begin{proof}
Choose an arbitrary coordinate system in which the first coordinate $u$ satisfies $\partial^a_u = \hat{N}^a$, then, up to a possible translation of $u$, there is a boundary defining function $\tilde{\varrho}$ with $\hat{N}^a\nabla_a\tilde{\varrho} =0$, such that:
\[ s= \frac{\tilde{\varrho}}{1+ u\sgn(\lambda)\tilde{\varrho}}. \]
Since $\varrho$ is a boundary defining function, $\tilde{\varrho} = e^\psi\varrho$ with $\mathcal{L}_{\hat{N}} \psi =0$, hence:
\[ s= \frac{\varrho}{e^{-\psi} + u\sgn(\lambda)\varrho}.\]
We obtain the desired coordinate system by setting:
\[ \tilde{u}= u - \sgn(\lambda)\varrho^{-1}(1-e^{-\psi}),\]
As $\tau = \lambda$ along the boundary, which implies that $\psi \to 0$ when $\varrho \to 0$, this is regular along the boundary.
\end{proof}
Whilst $\lambda$ is not an adapted scale, it turns out that the scaling of $n^a$ it provides nevertheless plays a distinguished role.
\begin{prop}\label{prop:distinguished_vector}
$\hat{N}^a$ is a conformal Killing vector, what is more for any choice of $\tau$ in Proposition~\ref{prop:coordinate_step1}, it is a Killing vector for $\bar{g}_{ab}=|\tau|^{-1}\zeta_{ab}$. 
\end{prop} 
\begin{proof}
From~\cite[Lemma 3.13 (3.18)]{RSTA20230042}\footnote{There are several typos in (3.18), it should be:
\[ \mathcal{L}_Ng^{ab}=-2N^{(a}\nu^2\lambda^{b)} + 2\nu\rho g^{ab}\]}, we have:
 \[\mathcal{L}_{\hat{N}} \hat{g}_{ab}= 2\hat{g}_{c(a}\hat{\nabla}_{b)}\hat{N^c}=-2s\sgn(\lambda)\hat{g}_{ab}.\]
 Since: $\bar{g}_{ab}=\frac{\varrho^2}{s^2}\hat{g}_{ab}$ we have:
\[\begin{aligned} \mathcal{L}_{\hat{N}}\bar{g}_{ab}&= 2\frac{\varrho}{s}\hat{N}^c\left(\frac{\nabla_c \varrho}{s}- \nabla_c s\frac{\varrho}{s^2}\right)\hat{g}_{ab}-2\frac{\varrho^2}{s}\sgn(\lambda)\hat{g}_{ab} \\ &=2\frac{\varrho}{s}\sgn(\lambda)\varrho\hat{g}_{ab}-2\frac{\varrho^2}{s}\sgn(\lambda)\hat{g}_{ab} =0. \end{aligned} \]
\end{proof}

\subsection{Normal form of the ambient metric near points of $\Sigma^\pm$}
We will now complete the study of the ambient structure in deriving a normal form for $\bar{g}_{ab}$ for some fixed scale $\tau$ satisfying the hypotheses of Prop.~\ref{prop:coordinate_step1} and use it to define a  distinguished local coordinate chart.
We begin with the following:
\begin{lemm}\label{lemm:metric_gbar}
\[\bar{g}_{ab}=-|\tau|^{-1}\sigma^2\sgn(\lambda)\nabla_a u \nabla_b u + 2\nabla_{(a} u\bar{\nabla}_{b)}\sigma|\tau|^{-\frac{1}{2}} + \sigma^{-2}(-\sgn(\lambda)+|\tau|^{-1}\nu^{-1})\bar{\nabla}_a\sigma\bar{\nabla}_b\sigma + \bar{h}_{ab}. \]
with:
\[\bar{h}_{ab}n^a= \bar{h}_{ab}t^b=0\]
and:
\[\mathcal{L}_{\hat{N}} \bar{h} =0.\]
\end{lemm}
\begin{proof}
Note that: 
\[\bar{g}_{ab}\hat{N}^b=|\tau|^{-1}|\lambda|\hat{\nabla}_a s= |\tau|^{-\frac{1}{2}}\bar{\nabla}_a\sigma -\sigma^2\tau^{-1}\nabla_a u, \]
plugging this into the right-hand side leads to:
\[ -\tau^{-1}\sigma^2\nabla_a u + |\tau|^{-\frac{1}{2}}\bar{\nabla}_a \sigma +\bar{h}_{ab}\hat{N}^b. \]
Comparing the two formulae we conclude that:
\[\bar{h}_{ab}\hat{N}^b=0.\]
Similarly, by definition:
\[\bar{g}_{ab}t^b= |\tau|^{-1}\sigma^{-2}\bar{\nabla}_a\sigma, \]
and plugging into the righthand side we have:
\[ -|\tau|^{-\frac{1}{2}}\nu \nabla_a u +|\tau|^{-\frac{1}{2}}\nu \nabla_a u +\nu\sgn(\lambda)\sigma^{-2}\bar{\nabla}_a\sigma + \sigma^{-2}(-\sgn(\lambda)\nu +|\tau|^{-1})\bar{\nabla}_a\sigma + \bar{h}_{ab}t^b.  \]
Hence, comparing we conclude that:
\[\bar{h}_{ab}t^b=0. \]
The fact that $\mathcal{L}_{\hat{N}}h=0$ follows from Proposition~\ref{prop:distinguished_vector}, and:
\[ \mathcal{L}_{\hat{N}} \varrho =0,\qquad \mathcal{L}_{\hat{N}}\nabla_a \varrho = 0, \qquad \mathcal{L}_{\hat{N}}\nabla_a u =0. \]
\end{proof}

We now check that the decomposition is coherent with the structure on the base, for this define:
\[ q_{ab}= \nu^{-1}|\tau|^{-1}\sigma^{-2} \bar{\nabla}_a\sigma\bar{\nabla}_b\sigma +\bar{h}_{ab}, \]
and recall that:
\[ g^{ab}=\sigma^2\zeta^{ab}=-lN^aN^b + \tilde{h}^{ab}= -\sgn(\lambda)\hat{N}^a\hat{N}^b + \tilde{h}^{ab}, \] 
with:
\[\mathcal{L}_{N} \tilde{h}^{ab} =0.\]
This tensor was used in~\cite{RSTA20230042} to induce the bilinear form $h^{ab}$ on $\mathcal{Q}$.

\begin{prop}
$\sigma^{-2}|\tau| q_{ab}$ is the pullback along the projection map $\pi$ of the induced metric on the quotient $\mathcal{Q}$.
\end{prop}
\begin{proof}
Let us calculate:
\[\begin{aligned} \sigma^{-2}|\tau|\tilde{h}^{ab}q_{bc}\\&=(\bar{g}^{ab} + \sgn(\lambda)\sigma^{-2}|\tau|\hat{N}^a\hat{N}^b)(\nu^{-1}|\tau|^{-1}\sigma^{-2}\bar{\nabla}_a\sigma\bar{\nabla}_b \sigma + \bar{h}_{ab})\\&=\bar{g}^{ab}(\bar{g}_{bc}+|\tau|^{-1}\sigma^2\sgn(\lambda)\nabla_b u \nabla_c u -2|\tau|^{-\frac{1}{2}}\nabla_{(b}u\bar{\nabla}_{c)}\sigma + \sgn(\lambda)\sigma^{-2}\bar{\nabla}_b\sigma\bar{\nabla}_c\sigma)\\&=\delta^{a}_c+ \sigma^2\sgn(\lambda)v^a\nabla_c u - |\tau|^{\frac{1}{2}}\sigma^2 t^a\nabla_c u + |\tau|\sgn(\lambda)t^a\bar{\nabla}_c\sigma - |\tau|^{\frac{1}{2}}v^a\bar{\nabla}_c\sigma.  \end{aligned}\] 
Where we have defined:
\[v^a=\zeta^{ab}\nabla_b u = (|\tau|^\frac{1}{2}\sigma^2 t^a -\hat{N}^a)\sigma^{-2}\sgn(\lambda).\]
Hence:
\[\begin{aligned} \sigma^{-2}|\tau|\tilde{h}^{ab}q_{bc}&=\delta^{a}_c -\hat{N}^a\nabla_c u + \hat{N}^a(\sigma^{-2}\sgn(\lambda)|\tau|^{\frac{1}{2}})\hat{\nabla}_a\sigma, \\&=\delta^a_c -\hat{N}^a(\nabla_c u -\sigma^{-2}\sgn(\lambda)|\tau|^\frac{1}{2}\hat{\nabla}_a\sigma),\\&=  \delta^{a}_c - N^a\frac{1}{2}\nabla_a l.\end{aligned}\]

Where we recall that: $l=\sigma^{-2}\lambda$ from which it follows that:
\[ \frac{1}{2}\nabla_a l = |\lambda|^{\frac{1}{2}}\sigma^{-1}(\nabla_a u - |\tau|^\frac{1}{2}\sigma^{-2}\hat{\nabla}_a\sigma \sgn(\lambda)).\]
\end{proof}

To conclude this section, we give the appropriate coordinate expression arising from the geometric decomposition of Lemma~\ref{lemm:metric_gbar}:
\begin{prop}\label{prop:final_coordinate}
Near any point of $\Sigma^\pm$, given a scale $\tau$ satisfying the hypotheses of Theorem~\ref{prop:coordinate_step1}, one can choose a coordinate system $(u,\varrho, y)$ such that:
\[ \frac{\partial^a}{\partial u} = \hat{N}^a, \qquad \frac{\partial^a}{\partial \varrho}=|\tau|^{\frac{1}{2}}\nu^{-1}\bar{t}^a, \qquad \sigma |\lambda|^{-\frac{1}{2}}=\frac{\varrho}{1+\sgn(\lambda)\varrho u}. \]
where $\bar{t}^a$ is defined by Lemma~\ref{lemm:properties_scale}.
In such a coordinate system the metric is of the form:
\[\begin{aligned} \bar{g}_{ab}&=-\varrho^2\sgn(\lambda)\nabla_a u \nabla_b u + 2\nabla_{(a} u \nabla_{b)}\varrho+ \varrho^{-2}(-\sgn(\lambda)+|\tau|^{-1}\nu^{-1})\nabla_a\varrho\nabla_b\varrho + \bar{h}_{ij}\nabla_a y^i\nabla_b y^j,\\ \bar{g}^{ab}&=(|\tau|\nu - \sgn(\lambda))\varrho^{-2}\hat{N}^a\hat{N}^b +2\sgn(\lambda)|\tau|\nu\hat{N}^{(a}{\partial^{b)}_\varrho}+\varrho^2|\tau|\nu{\partial^a_\varrho}\partial^b_\varrho+\bar{h}^{ij}\partial^a_{y^i}\partial^b_{y^j}.   \end{aligned} \]
\end{prop}
\begin{rema}
We remark that in the presence of the infamous mass aspect, we would only recover a $b$-metric at this stage.
\end{rema}
\begin{proof}
Let $(u,x)$ be an arbitrary coordinate system near $p\in \Sigma^\pm$ given by Proposition~\ref{prop:coordinate_step1}.
Since $\hat{N}^a\nabla_a\varrho =0$ and $\nabla_a\varrho \neq 0$ on $\Sigma^\pm$, permuting the variables if necessary one can assume that $\frac{\partial \rho}{\partial x_0}\neq 0$.
Consider then the coordinate transformation:
\[ c: (u,x_0,y) \mapsto (u,\varrho,y),\]
where $ x= (x_0,y)$. The Jacobian matrix of this transformation is given by:
\[ \begin{pmatrix} 1 & 0 & 0\\ 0 & \frac{\partial \varrho}{\partial x_0} & \frac{\partial \varrho}{\partial y} \\ 0 & 0 &  \delta^i_j, \end{pmatrix}\]
which is invertible. Hence $(u,\varrho, y)$ can be used as a coordinate chart near $p$.
By Lemma~\ref{lemm:metric_gbar}, one can see that $\{\hat{N}^a, |\tau|^\frac{1}{2}\nu^{-1}t^a\}^\perp=\textrm{span}(\frac{\partial^a}{\partial y^i})$. Using again Lemma~\ref{lemm:metric_gbar}, and decomposing $\frac{\partial^a}{\partial \varrho}$ with respect to the above orthogonal decomposition away from $\Sigma$ we can see that:
\[\begin{aligned} \frac{\partial^a}{\partial \varrho} &= -|\tau|\sigma^{-2}\sgn(\lambda)\hat{N}^a + |\tau|^{\frac{1}{2}}\nu^{-1}t^b + \sigma^k(u,\varrho,y)\frac{\partial^a}{\partial y^k},\\&=|\tau|^\frac{1}{2}\nu^{-1}\bar{t}^a + \sigma^k(u,\varrho,y)\frac{\partial^a}{\partial y^k}. \end{aligned}\]
The final expression extends to $\Sigma$ by continuity. Recall that:
\[\bar{\nabla}_c n^b= 0, \qquad \bar{\nabla}_n \nu =0,\]
therefore:
\[\begin{aligned} \mathcal{L}_{\hat{N}}(|\tau|^\frac{1}{2}\nu^{-1}\bar{t}^a)&= |\lambda|^\frac{1}{2}|\tau|^\frac{1}{2}\nu^{-1}(n^a\bar{\nabla}_a \bar{t}^b  - \frac{1}{2}\lambda^{-1}\bar{t}^a\bar{\nabla}_a\lambda n^b).\end{aligned}\]
Now:
\[ \sigma \bar{\nabla}_n\bar{t}^b = \nabla_n(-\lambda^b + \sigma^{-1}(1-\sgn(\lambda)|\tau|^\frac{1}{2}|\lambda|^\frac{1}{2}\nu)n^b), \]
with:
\[ \nabla_n\bar{\lambda}^b= n^a\nabla_a(H^{AB}Y_AZ_B^b)=-\nu n^b,\]
so that:
\[\sigma \bar{\nabla}_n \bar{t}^b = \nu n^b  -\frac{|\tau|^\frac{1}{2}}{|\lambda|^\frac{1}{2}}\nu n^b=\nu(1-\frac{|\tau|^\frac{1}{2}}{|\lambda|^\frac{1}{2}})n^b. \]
Additionally:
\[ \begin{aligned} t^a \bar{\nabla}_a \lambda &= \sigma^{-2} \bar{\nabla}_b\sigma \zeta^{ab}\bar{\nabla}_a\lambda\\&=\sigma^{-2}W_b^CD_C\sigma H^{AB}Z^a_AZ^b_BW_a^DD_D\lambda\\&= \sigma^{-2}(\delta^C_B - X^CY_B)(\delta^D_A-X^DY_A)H^{AB}D_D\lambda D_C\sigma\\&=\sigma^{-2}(2\sigma +2\sigma\lambda\nu-2\sigma)=2\sigma^{-1}\nu\lambda.\end{aligned}\]
From which it follows that:
\[\bar{t}^a\bar{\nabla}_a \lambda = (t^a-\sigma^{-2}\sgn(\lambda)|\tau|^\frac{1}{2}|\lambda|^\frac{1}{2}\nu n^a)\bar{\nabla}_a\lambda=2\sigma^{-1}\nu\lambda(1-\frac{|\tau|^\frac{1}{2}}{|\lambda|^\frac{1}{2}}).  \]
Plugging this into the above leads to the conclusion that:
\[ \mathcal{L}_{\hat{N}}(|\tau|^\frac{1}{2}\nu^{-1}\bar{t}^a) =0. \]
Since,
\[ \mathcal{L}_{\hat{N}} \partial^a_\varrho =0, \]
we conclude that $\sigma^k(u,\varrho,y)\equiv \sigma^k(\varrho,y)$ and then the transformation, 
\[ (u,\varrho,y) \mapsto (u, \varrho, \tilde{y}^k= y^k+ \int^\varrho \sigma^k(s,y)\textrm{d}s), \]
yields a coordinate system with the desired property.
\end{proof}

\section{The projectively invariant Klein-Gordon system}\label{sec:KG}
\subsection{Klein-Gordon equation}
Having now studied the local geometry near the extended boundary $\Sigma^\pm$,  we turn our attention to the heart of our paper: the description of a framework for the study of the Klein-Gordon equation on the base $\mathcal{Q}_{\mathcal{O}}$. The general philosophy follows that of Penrose's conformal method by isolating a projectively invariant system (with respect to the ambient projective structure) that imposes the Klein-Gordon equation when expressed in the appropriate \enquote{physical} scale. On the base manifold $\Q$, this scale is by definition $\bar{\sigma}$ which pulls back to the density $\sigma$ on the ambient space $M$. 

For definiteness, we will consider the physical manifold to be given by $\Q_{\mathcal{O}}^+=\{\bar{\sigma} >0\}$, and $\Q^+=\{\bar{\sigma}\geq 0\}$ its projective compactification. This is naturally covered by $\mathcal{O}^+=\{\sigma >0\}$, and $M^+=\{\sigma\geq 0\}$ respectively. On $\mathcal{O}_+$, the scale $\sigma$ determines a unique distinguished connection that we denote by $\nabla$ and that enjoys the following properties:
\begin{prop}[{\cite[Prop 4.4]{RSTA20230042}}]
\begin{equation}\label{eq:essential_properties_scale_sigma} \nabla_c n^a = 0, \qquad \nabla_c \zeta^{ab}=-2\sigma^{-1}\delta_c^{(a}n^{b)}, \quad P_{ab}=0. \end{equation}
\end{prop}
The following \emph{scale-dependent} definitions will be useful in our discussion:
\begin{equation} g^{ab}=\sigma^2\zeta^{ab} \qquad N^a=\sigma n^a=\sigma Z^a_AI^A, \qquad l=\sigma^{-2}\lambda \tag{{\cite[(4.1)]{RSTA20230042}}}, \end{equation}
We recall once more that:
\begin{equation} \label{eq:sigma_scale_metric_decomposition} g^{ab}=\sigma^2\zeta^{ab}=-lN^aN^b + \tilde{h}^{ab}, \qquad \mathcal{L}_{N}\tilde{h}^{ab}=0 \qquad \tilde{h}^{ab}\nabla_b l =0, \end{equation}
and that $\tilde{h}^{ab}$ pushes forward to the metric on $\Q^+_{\mathcal{O}}$.  However, we emphasise that $\nabla$ is \emph{not} the Levi-Civita connection of $g$, which we will denote by $\tilde{\nabla}$. They are however related by (see ~\cite[Eq.(4.12)]{RSTA20230042}):
\begin{equation}\label{eq:LinkScaleSigmaLeviCivita_vectors} \tilde{\nabla}_c \xi^b =\nabla_c \xi^b + N^bg_{cd}\xi^d , \end{equation}
we note that this implies respectively for a form $\omega_b$ and a weight $w$ density $\tau$:
\begin{equation}\label{eq:LinkScaleSigmaLeviCivita_forms_densities} \begin{aligned} \tilde{\nabla}_c \omega_b= \nabla_c \omega_b -N^a\omega_ag_{cb} \\ \tilde{\nabla}_c \tau = \nabla_c \tau + \frac{w}{n+2}g_{cd}N^d\tau.\end{aligned}\end{equation}

Finally, written in the scale $\nabla$ the ambient tractor metric on $M^+$ takes the form:
\begin{equation}\label{eq:ambient_tractor_metric_scale_sigma} H^{AB}= \zeta^{ab}W_a^AW_b^B  +2\sigma^{-1}n^bX^{(A}W_b^{B)}, \end{equation}
and we recall once more that:
\[ x \in \mathcal{O}^+ \mapsto (\frac{1}{2}{l(x)}, \pi(x)) \in \R \times \mathcal{Q}^+_\mathcal{O}, \]
is a distinguished global trivialisation and $N$ is the fundamental vector field $1^*$. 

It turns out that the projectively invariant operator encoding a Klein-Gordon equation is the natural one and that the main subtlety lies in restricting our attention to an appropriate subspace of densities, this is the content of the key result in this section:
\begin{theo}\label{thm:proj_KG}
Let $m>0$ be fixed and $\phi$ be a density of weight $w=-\frac{n-1}{2}$.  Assume that $\phi$ satisfies the projectively invariant system:
\begin{equation}\label{eq:proj_KG}\begin{cases}H^{AB}D_AD_B\phi =0,\\ (D_n -\frac{im}{\sqrt{|\lambda|}})\frac{\sqrt{|\lambda|}}{2im}(D_n +\frac{im}{\sqrt{|\lambda|}})\phi=0, \end{cases}\end{equation}
Then, on $\mathcal{O}^+=\mathcal{O}\cap\{\sigma >0\}$, setting $\phi= \varphi \sigma^{-\frac{n-1}{2}}$ and \[ \Box^h = \tilde{h}^{ab}\left(\delta^c_a - \frac{\sigma^{-1}}{2}\nabla_a \lambda n^c \right)\tilde{\nabla}_c \left(\delta^d_b - \frac{\sigma^{-1}}{2}\nabla_b \lambda n^d\right)\tilde{\nabla}_d=\tilde{h}^{ab}\tilde{\nabla}_a\tilde{\nabla}_b\] we have:
\[(\Box^h +\sgn(\lambda) m^2)\varphi=0.\]
\end{theo}

Let us comment first on the definition of the operator $\Box^h$. It is written on the ambient space, but is meant to be understood to act like the intrinsic Laplacian operator on the base manifold $\mathcal{Q}^+_{\mathcal{O}}$. We emphasise that this intuition is correct and justified by the explicit appearance in the above of the horizontal projection (using the canonical principle connection): 
\[ \pi^c_{\,\,a}=\delta^c_a - \frac{\sigma^{-1}}{2}\nabla_a \lambda n^c\,,\]
which shows that the operator is acting in the horizontal directions. This point will again be consolidated by Proposition~\ref{prop_commutation_of_y_Dn}, which is a projective version of the statement that $y$ commutes with the \enquote{Lie derivative in the direction $n^a$}.

It should be noted that the field $\varphi$ in Theorem~\ref{thm:proj_KG} is a field on the ambient manifold, and it's Lie derivative in the fibre direction $N^a$ does \emph{not} vanish; it has a particular dependence that we will study in the Section~\ref{sec:m_density_space}; we will see that solutions of the Klein-Gordon equation on the base are obtained by pulling back $\varphi$ along constant sections of $l=\sigma^{-2}\lambda$.
We prove the theorem in several steps, first we justify that the two expressions for $\Box^h$ coincide.
\begin{lemm}
 \[\tilde{h}^{ab}\left(\delta^c_a - \frac{\sigma^{-1}}{2}\nabla_a \lambda n^c \right)\tilde{\nabla}_c \left(\delta^d_b - \frac{\sigma^{-1}}{2}\nabla_b \lambda n^d\right)\tilde{\nabla}_d=\tilde{h}^{ab}\tilde{\nabla}_a\tilde{\nabla}_b\] 
\end{lemm}
\begin{proof}
By definition of $\tilde{h}^{ab}$ we have $\tilde{h}^{ab}\nabla_b\lambda =0$, hence the first projection operator is already unnecessary, i.e.
\[ \tilde{h}^{ab}\left(\delta^c_a - \frac{\sigma^{-1}}{2}\nabla_a \lambda n^c \right)\tilde{\nabla}_c \left(\delta^d_b - \frac{\sigma^{-1}}{2}\nabla_b \lambda n^d\right)\tilde{\nabla}_d= \tilde{h}^{cb}\tilde{\nabla}_c \left(\delta^d_b - \frac{\sigma^{-1}}{2}\nabla_b \lambda n^d\right)\tilde{\nabla}_d.\]
We expand out the remaining part to find:
\[ \tilde{h}^{cd}\tilde{\nabla}_c\tilde{\nabla}_d - \tilde{h}^{cb}\tilde{\nabla}_c\left(\frac{\sigma^{-1}}{2}\nabla_b\lambda n^d\right)\tilde{\nabla}_d -\tilde{h}^{cb}\left(\frac{\sigma^{-1}}{2}\nabla_b\lambda n^d\right)\tilde{\nabla}_c\tilde{\nabla}_d.\]
The final term vanishes as before and hence we wish to show that the middle term vanishes.
Using Eqs.~\eqref{eq:LinkScaleSigmaLeviCivita_vectors} and \eqref{eq:LinkScaleSigmaLeviCivita_forms_densities}, we see that:
\[\begin{aligned} \tilde{h}^{cb}\tilde{\nabla}_c\left(\frac{\sigma^{-1}}{2}\nabla_b\lambda n^d\right)&=\tilde{h}^{cb}\left(\nabla_c\left(\frac{\sigma^{-1}}{2}\nabla_b\lambda n^d\right) + \frac{\sigma^{-1}}{2}\nabla_b\lambda N^d g_{ca}n^a - \frac{\sigma^{-1}}{2}N^a\nabla_a\lambda n^dg_{cb} \right)\\&= \tilde{h}^{cb}\nabla_c\left(\frac{\sigma^{-1}}{2}\nabla_b\lambda n^d\right)-\sigma n^d \tilde{h}^{cb}g_{cb},\end{aligned}\]
where we have used that: $N^a\nabla_a \lambda = 2 \sigma^2.$
Since $\sigma$ and $n^a$ are both parallel for the connection $\nabla$ and by Eqs.~\eqref{eq:CsqParallel} :
\[ \frac{1}{2}\nabla_c \nabla_b \lambda = \zeta_{cb}-\sigma^{-2}\lambda n_cn_b,\]
it follows that:
\[ \tilde{h}^{cb}\tilde{\nabla}_c\left(\frac{\sigma^{-1}}{2}\nabla_b\lambda n^d\right)=0, \]
which concludes the proof.
\end{proof}
Next we compute the projectively invariant operator $H^{AB}D_AD_B$:
\begin{prop} \label{eq:laplace_in_sigma_scale}
Let $\phi\in\mathcal{E}(w)$ be projective density and let $\varphi:= 	\phi \sigma^{-w}$ be its expression in the scale $\sigma$. Then
	\begin{equation}
		H^{AB}D_AD_B\phi=\sigma^{w-2} \Big( \tilde{\Box} \varphi +\left(2w + n -1\right) N^a\nabla_a \varphi\Big).
	\end{equation}
In particular, in the scale $\sigma$ and when acting on densities of weight $w=\frac{1-n}{2}$, $H^{AB}D_AD_B\phi$ coincides with the Laplacian of $g$.
\end{prop}
\begin{proof}
A standard computation shows that in an arbitrary scale $\hat{\nabla}$:
\[ D_AD_B\phi= w(w-1)\phi \hat{Y}_A\hat{Y}_B + 2(w-1)\nabla_a\phi Z^a_{(A}\hat{Y}_{B)} + (\hat{\nabla}_a \hat{\nabla}_b \phi + w \phi \hat{P}_{ab})Z^a_AZ^b_B.\]
	Taking into account the properties~\eqref{eq:essential_properties_scale_sigma} of the scale $\nabla$ we have:
	\begin{align*}
		H^{AB}D_AD_B\phi&= \zeta^{ab}\nabla_a\nabla_b\phi +2(w-1)\nabla_a \phi \sigma^{-1}n^a,\\
		& = \sigma^{w-2}\Big(g^{ab}\nabla_a\nabla_b\varphi +2(w-1)N^a\nabla_a \varphi \Big)
	\end{align*}
	Now, from~\eqref{eq:LinkScaleSigmaLeviCivita_forms_densities} it follows that:
	\[\tilde{\Box}\varphi = g^{ab}\nabla_a\nabla_b\varphi -(n+1)N^a\nabla_a \varphi. \]
	Putting everything together yields the result.
\end{proof}

The next computational lemma explores the condition $\phi \in \mathcal{E}_m(w)$:
\begin{lemm}\label{lemm:csq_condition2}
\[\begin{aligned} \sigma(D_n -\frac{im}{\sqrt{|\lambda|}})\sqrt{|\lambda|}(D_n +\frac{im}{\sqrt{|\lambda|}})\phi=|l|^\frac{1}{2}N^a\nabla_a(N^a\nabla_a\phi)+\frac{\sgn(l)}{|l|^\frac{1}{2}}N^a\nabla_a\phi +\frac{m^2}{|l|^\frac{1}{2}}\phi.\end{aligned}\]
\end{lemm}
\begin{proof}
\[\begin{aligned} \sigma(D_n -\frac{im}{\sqrt{|\lambda|}})\sqrt{|\lambda|}(D_n +\frac{im}{\sqrt{|\lambda|}})\phi &
=(N^a\nabla_a - \frac{im}{|l|^\frac{1}{2}})(|l|^\frac{1}{2}N^a\nabla_a\phi +im\phi)\\&= N^a\nabla_a(|l|^\frac{1}{2}N^a\nabla_a\phi) + imN^a\nabla_a\phi -imN^a\nabla_a\phi +\frac{m^2}{|l|^\frac{1}{2}}\\&=|l|^\frac{1}{2}N^a\nabla_a(N^a\nabla_a\phi)+\frac{\sgn(l)}{|l|^\frac{1}{2}}N^a\nabla_a\phi +\frac{m^2}{|l|^\frac{1}{2}}\phi.\end{aligned}\]
\end{proof}
We can now prove Theorem~\ref{thm:proj_KG}:
\begin{proof}[Proof of Theorem~\ref{thm:proj_KG}]
Starting from the result of Proposition~\ref{eq:laplace_in_sigma_scale}, we already know that when $\phi \in \mathcal{E}(-\frac{n-1}{2})$, and if $\varphi=\phi\sigma^\frac{n-1}{2}$, then:
\[ H^{AB}D_AD_B \phi= \sigma^{-\frac{n+3}{2}}\tilde{\Box}\varphi. \]
Now equating two different expressions for $\Phi_{AB}$ using ~\eqref{eq:sigma_scale_metric_decomposition} and~\eqref{eq:ambient_tractor_metric_scale_sigma}:
\[\begin{aligned} \tilde{\Box}\varphi= g^{ab}\tilde{\nabla}_a\tilde{\nabla}_b \varphi&= -lN^aN^b\tilde{\nabla}_a\tilde{\nabla}_b\varphi + \Box^h\varphi \\&= -lN^a\tilde{\nabla}_{a}(N^b\tilde{\nabla}_b\varphi) + \Box^h \varphi - N^b\tilde{\nabla}_b\varphi,\end{aligned}\]
where we have used the fact that:
\[\tilde{\nabla}_a N^b = N^bg_{ac}N^c.\]

Now multiplying the result of Lemma~\ref{lemm:csq_condition2}, by $|l|^\frac{1}{2}\sgn(l)$ and using this in the above we arrive at:
\[ H^{AB}D_AD_B \phi = \sigma^{-\frac{n+3}{2}}(\Box^h+ \sgn(\lambda) m^2)\varphi.\]
\end{proof}

\subsection{The mass-shell density space $\mathcal{E}_m(w)$}\label{sec:m_density_space}
In this section, we shall focus our attention on the structure of the density subspace in which we look for our solutions: \[\mathcal{E}_m(w)=\ker~(D_n -\frac{im}{\sqrt{|\lambda|}})\frac{\sqrt{|\lambda|}}{2im}(D_n +\frac{im}{\sqrt{|\lambda|}}).\]
One can distinguish two subspaces of $\mathcal{E}_m(w)$, namely:
\[ \mathcal{E}^\pm_m(w)= \ker\left(D_n \mp \frac{im}{\sqrt{|\lambda|}}\right).\] 
It turns out they generate they entire space:
\begin{prop}\label{prop:decomposition_champs}
\[ \mathcal{E}_m(w)= \mathcal{E}_m^+(w)\oplus \mathcal{E}_m^-(w). \]
\end{prop}
\begin{proof}
Note that:
\[ \begin{gathered} [D_n -\frac{im}{|\lambda|^\frac{1}{2}}, \frac{|\lambda|^\frac{1}{2}}{2im}]=\frac{\sgn(\lambda)\sigma}{2im |\lambda|^\frac{1}{2}}, \quad
[D_n +\frac{im}{|\lambda|^\frac{1}{2}},D_n - \frac{im}{|\lambda|^\frac{1}{2}}]=-2im[D_n,\frac{1}{|\lambda|^\frac{1}{2}}]=\frac{2im\sigma\sgn(\lambda)}{|\lambda|^{\frac{3}{2}}} \end{gathered}\]
Hence:
\[ \begin{aligned} (D_n - \frac{im}{\sqrt{|\lambda|}})\frac{\sqrt{|\lambda|}}{2im}(D_n + \frac{im}{\sqrt{|\lambda|}})&= \frac{\sgn(\lambda)\sigma}{2im|\lambda|^\frac{1}{2}}(D_n+\frac{im}{\sqrt{|\lambda|}})+\frac{\sqrt{|\lambda|}}{2im}(D_n-\frac{im}{\sqrt{|\lambda|}})(D_n+\frac{im}{\sqrt{|\lambda|}}),\\&=\frac{\sgn(\lambda)\sigma}{2im|\lambda|^\frac{1}{2}}(D_n-\frac{im}{\sqrt{|\lambda|}})+\frac{\sqrt{|\lambda|}}{2im}(D_n+\frac{im}{\sqrt{|\lambda|}})(D_n-\frac{im}{\sqrt{|\lambda|}}),\\&=(D_n + \frac{im}{\sqrt{|\lambda|}})\frac{\sqrt{|\lambda|}}{2im}(D_n - \frac{im}{\sqrt{|\lambda|}}). \end{aligned}\]
Therefore: 
\[ \mathcal{E}^\pm_m(w)\subset \mathcal{E}_m(w). \]
%
Furthermore, if we set: \[p_\pm = \pm\frac{|\lambda|^\frac{1}{2}}{2im}(D_n\pm \frac{im}{\sqrt{|\lambda|}}), \]
then \[ p_+ + p_- = 1,\]
and on $\mathcal{E}_m(w)$: 
\[ p_+p_-=p_-p_+=0,\]
which proves that:
\[\mathcal{E}_m(w)=\ker p_+ \oplus \ker p_- = \mathcal{E}^-_m(w)\oplus \mathcal{E}^+_m(w).\] 
\end{proof}
In fact, we have an important stability result:
\begin{prop}\label{prop:stability_decomposition}When $w=-\frac{n-1}{2}$,the direct sum decomposition of Proposition~\ref{prop:decomposition_champs} is stable under the projectively invariant operator, $H^{AB}D_AD_B$, in other words:
\[ [D_n \pm \frac{im}{|\lambda|^\frac{1}{2}},H^{AB}D_AD_B]=0.\] 
\end{prop}
This result is an immediate consequence of Lemma~\ref{prop_commutation_lambda_and_y} and Proposition~\ref{prop_commutation_of_y_Dn}.  These results are of  independent interest that will be discussed in Section~\ref{sec:kg_boundary_calculus}.

The significance of Proposition~\ref{prop:stability_decomposition} lies in the fact that:
\[ H^{AB}D_AD_B \phi = 0 \Leftrightarrow\begin{cases} H^{AB}D_AD_B \phi_\pm =0,\\ \phi_\pm = p_\pm\phi,\end{cases}\]
meaning we can solve the ambient Klein-Gordon system~\eqref{eq:proj_KG} by working on the components of the direct sum decomposition~\eqref{prop:decomposition_champs}  separately. This amounts to studying:
\begin{equation}\label{eq:proj_KGpm}\begin{cases}H^{AB}D_AD_B\phi_\pm =0,\\ (D_n \mp \frac{im}{\sqrt{|\lambda|}})\phi_\pm=0. \end{cases} \tag{\ref{eq:proj_KG}$_\pm$}\end{equation}
Working in the scale $\sigma$, the second condition in~\eqref{eq:proj_KGpm} simply states that the functions: $\varphi_\pm =\phi_\pm\sigma^{\frac{n-1}{2}}$ satisfy:
\[ \mathcal{L}_N \varphi_\pm = \pm\frac{ im}{|\lambda|^{\frac{1}{2}}}.\]
This is immediately solved and we deduce that:
\[ \phi = (a_+e^{im\sgn(\lambda)\sqrt{|l|}} + a_-e^{-im\sgn(\lambda)\sqrt{|l|}})\sigma^{-\frac{n-1}{2}}, \]
where both $a_+$ and $a_-$ satisfy $\mathcal{L}_{N} a_\pm =0$ and by Theorem~\ref{thm:proj_KG} are the pullbacks of solutions to the Klein-Gordon equation on the base $\mathcal{Q}_{\mathcal{O}}.$

\subsection{Asymptotics of solutions of the Klein-Gordon equation}
In this section we aim to study the behaviour of solutions of the ambient projective Klein Gordon system~\eqref{eq:proj_KG} near the boundary; for this, we will exploit the results of Section~\ref{sec:ambient_geometry}. 
Throughout, we fix an adapted scale $\tau$ as in Proposition~\ref{prop:coordinate_step1} and use the associated local coordinate chart $(u,\varrho,y)$ from Proposition~\ref{prop:final_coordinate}. We denote by $\bar{\nabla}$ the connection in the projective class determined by $\tau$ and reintroduce the tensor fields from Eq.~\eqref{eq:tensors_scale_tau}:
\[\bar{g}_{ab}=|\tau|^{-1}\zeta_{ab},\qquad \nu = \sigma^{-2}\zeta^{ab}\bar{\nabla}_a\sigma\bar{\nabla}_b\sigma,\qquad t^b=\sigma^{-2}\zeta^{ab}\bar{\nabla}_a\sigma. \tag{\ref{eq:tensors_scale_tau}} \]

In the projective splitting given by the scale $\bar{\nabla}$, the tractor metric takes the following form
\[ H^{AB}= \zeta^{ab}\bar{W}^A_a\bar{W}^B_b + 2\lambda^bX^{(A}\bar{W}^{B)}_b + \nu X^AX^B,\]
with:
\[ \lambda^a=H^{AB}Y_BZ_A^a= \sigma^{-1}(1-\sgn(\lambda)|\tau|^\frac{1}{2}|\lambda|^\frac{1}{2}\nu)n^b - \sigma \bar{t}^b.\]
Similarly to the scale $\sigma$, the connection $\bar{\nabla}$ determined by $\tau$ in the projective class does \emph{not} coincide with the Levi-Civita connection $\nabla^{\bar{g}}$ of the metric $\bar{g}$, but they are closely related by the relation:
\[\begin{aligned} \nabla^{\bar{g}}_a\xi^b &= \bar{\nabla}_a\xi^b +\lambda^b \zeta_{ac}\xi^c.
 \end{aligned}\]
Using this one can derive the form of the operator $H^{AB}D_AD_B$ in such a scale:
\begin{prop}\label{eq:KGinBoundaryScale}
Let $\phi$ be an arbitrary section of $\mathcal{E}(-\frac{n-1}{2})$ then:
\[ H^{AB}D_AD_B \phi = |\tau|^{-\frac{n+3}{2}}\left(\Box^{\bar{g}}\psi -\frac{n-1}{2}\bar{g}^{ab}\bar{P}_{ab}\psi+\frac{(n-1)(n+1)}{4}\nu|\tau|\psi \right), \]
where $\psi=\phi|\tau|^{\frac{n-1}{2}}$.
\end{prop}
We will not develop further the local expression of this operator in this paper, as it is not required for our current purposes. By Prop.~\ref{prop:stability_decomposition}, we can of course restrict our attention to the subspaces $\mathcal{E}^+_m(-\frac{n-1}{2})$. Since $D_n\tau =0$, the second condition in Eq.~\eqref{eq:proj_KGpm} can be written in terms of the scalar function \[ \psi_\pm=\phi_\pm |\tau|^\frac{n-1}{2}, \]
as:
\[ \mathcal{L}_{\hat{N}}\psi_\pm= \pm im \psi_\pm, \]
since locally $\hat{N}= \partial_u,$ this is immediately solved and we find that:
\[ \psi_\pm = b_{\pm}e^{\pm im u}, \]
where the functions $b_{\pm}$ satisfy:
\[ \mathcal{L}_{\hat{N}}b_{\pm}= 0.\]

Restricting to $\psi_+$ for readability ($\psi_-$ is completely analogous), we can equate the two expressions for $\phi_+$ in the interior and boundary scales $\sigma$ and $\tau$ respectively:
\[ a_+e^{im\sgn(\lambda)\sqrt{|l|}}\sigma^{-\frac{n-1}{2}} =\phi_+= b_+e^{imu}|\tau|^{-\frac{n-1}{2}}.\]
Since $\sigma= \varrho |\tau|^\frac{1}{2}$, we arrive at the following equality in a neighbourhood of the boundary:
\[ a_+e^{im \sgn(\lambda)\sqrt{|l|}}\varrho^{-\frac{n-1}{2}}=b_+e^{imu} \Leftrightarrow b_+ =a_+\varrho^{-\frac{n-1}{2}}e^{im(\sgn(\lambda)\sqrt{|l|}-u)}. \]

Now, if we assume that the field $\phi_+$ extends continuously to the boundary, then $b_+$ must have a finite limit along $\varrho \to 0$. Since:
\[ \sgn(\lambda)\sqrt{|l|}-u= \frac{\sgn(\lambda)}{s} -u= (\frac{\sgn(\lambda)+u\varrho}{\varrho} -u)=\frac{\sgn(\lambda)}{\varrho};\]
this imposes that:
\[ \lim_{\varrho \to 0} \varrho^{\frac{1-n}{2}}e^{\sgn(\lambda)\frac{im}{\varrho}}a_+.\]
is finite. Summarising:
\begin{theo}\label{thm_asymptotics}
Let $\phi \in \mathcal{E}_m(-\frac{n-1}{2})$ solves $D_AD_B\phi=0$ and assume that the elements of its decomposition $\phi_+$ and $\phi_-$ have a finite limit along the boundary, then for any choice of $l_0\neq 0$ the solution $\varphi = a_+e^{im\sgn(\lambda)\sqrt{|l_0|}} + a_-e^{-im\sgn(\lambda)\sqrt{|l_0|}}$ of the Klein-Gordon equation on $\mathcal{Q}_{\mathcal{O}}$ 
 obtained by restricting $\phi$ to the section $\{l=l_0\}$, must be such that:
 \[ \lim_{\varrho\to 0} a_\pm\varrho^{\frac{1-n}{2}}e^{\pm\sgn(\lambda)\frac{im}{\varrho}}, \]
 are finite for any choice of boundary defining function $\varrho$ obtained from an adapted density $\tau$.
\end{theo}
These asymptotics correspond precisely with the classical result due to Hörmander \cite{hormanderRemarksKleinGordonEquation1987} (see: Introduction) for the \enquote{positive and negative frequency parts} of Klein-Gordon equation in Minkowski space time. The statement here suggests a partial peeling type result, in the spirit of that of Penrose~\cite{penrose_zero_1965} (see also~\cite{mason_nicolas_2009}) relating Hörmander's asymptotic behaviour to regularity assumptions of a field near the projective boundary, in the present case, in a higher-dimensional ambient space. Additionally, whereas Hörmander's decomposition follows from the Fourier transform of the field, here it arises structurally in the \emph{ambient}  space projective geometry. Nevertheless, it is only recovered asymptotically for the field on the base. To show this, we first we discuss how to \enquote{project out} the coefficients $a_\pm$; we begin with:
 \begin{lemm}\label{lemm:recover_decomp}
 Fix $l_0\in \R$ and define: \[ \bar{\varphi} = a_+e^{im\sgn(\lambda)\sqrt{|l_0|}} + a_-e^{-im\sgn(\lambda)\sqrt{|l_0|}}=\phi\sigma^{\frac{n-1}{2}}|_{l=l_0},\]
 where $\phi$ solves Eq.~\eqref{eq:proj_KG} then:
 \[ a_+e^{im\sgn(\lambda)\sqrt{|l_0|}} = \frac{1}{2}\left(\bar{\varphi} + \frac{\sgn(\lambda)\varrho^2}{im}\left(\frac{\partial \bar{\varphi}}{\partial \varrho}-\left.\frac{\partial (\phi\sigma^\frac{n-1}{2})}{\partial \varrho}\right|_{l=l_0} \right)\right).\]
  \end{lemm}
 \begin{proof}
 Recall the definition of $\phi_+ = p_+\phi=\frac{|\lambda|^\frac{1}{2}}{2im}(D_n +\frac{im}{\sqrt{|\lambda|}})\phi $. Since $D_n\sigma =0$, it follows that:
 \[ \phi_+\sigma^\frac{n-1}{2}= \frac{1}{2im}\partial_u(\phi\sigma^{\frac{n-1}{2}}) + \frac{1}{2}(\phi\sigma^{\frac{n-1}{2}}).  \] 
 Now: 
 \[ l=l_0 \Rightarrow \varrho u = -\sgn(\lambda) +\sgn(\lambda)\varrho\sqrt{|l_0|},\]
and the result follows from the chain rule and pulling back along the chosen section.
 \end{proof}
 
 We emphasise a major pitfall of this lemma: the final result still depends on the transverse derivative of the \emph{full} ambient space field $\phi\sigma^\frac{n-1}{2}$. However, this dependence vanishes in the limit $\varrho \to 0$, under the additional assumption that $b_+=e^{-imu}\phi_+|\tau|^{\frac{n-1}{2}}$ is at least $\mathcal{C}^1$ up to the boundary. Indeed, writing:
\[\phi_+= b_+e^{imu}|\tau|^{-\frac{n-1}{2}}\]
so that:
\[\phi_+\sigma^{\frac{n-1}{2}}=\varrho^{\frac{n-1}{2}} b_+ e^{imu} \]
taking the derivative we therefore obtain:
\[ \partial_\varrho \left(\phi_+\sigma^{\frac{n-1}{2}}\right)=\varrho^{\frac{n-3}{2}}\left(\frac{n-1}{2}b_++  \varrho\partial_\varrho b_+ \right)e^{imu}= O(\varrho^{\frac{n-3}{2}}), \]
from which it follows that:
\[ \frac{\partial (\phi_+\sigma^\frac{n-1}{2})}{\partial \varrho}\varrho^\frac{1-n}{2}\varrho^2e^{\sgn(\lambda)\frac{im}{\varrho}}=O(\varrho). \]
Therefore, overall:
\begin{theo}\label{theo:asymptotics2}
Let $\phi = \phi_+ + \phi_- \in \mathcal{E}_m(-\frac{n-1}{2})$ solve the Klein-Gordon system~\eqref{eq:proj_KG} and suppose that:
\[ b_+=\phi_+|\tau|^{\frac{n-1}{2}}, \]
is $\mathcal{C}^1$ up to the boundary. Then for any $l_0 \in \R$, the pullback $\bar{\varphi}$ of $\phi\sigma^{\frac{n-1}{2}}$ along the section $l=l_0$ is a solution of the Klein-Gordon equation and:
\begin{equation}\begin{aligned}\lim_{\varrho\to 0}\frac{1}{2} \left(\bar{\varphi} + \frac{\varrho^2}{im}\frac{\partial\bar{\varphi}}{\partial \varrho} \right)\varrho^\frac{1-n}{2}e^{\sgn(\lambda)\frac{im}{\varrho}}&=\lim_{\varrho \to 0} \frac{1}{2}\left(\bar{\varphi} + \frac{\nu^{-1}|\tau|^{-1}}{im}q^{ab}\nabla_b \varrho \nabla_a\bar{\varphi} \right)\varrho^{\frac{1-n}{2}}e^{\sgn(\lambda)\frac{im}{\varrho}} \\&=\lim_{\varrho \to 0} b_+ e^{-im\sgn(\lambda)\sqrt{|l_0|}}. 
\end{aligned}\end{equation}
An analogous statement holds for $\phi_-$.
\end{theo}

In this paper, we restrict ourselves to extracting the asymptotics of the solutions to the Klein-Gordon equation that are obtained through this higher-dimensional framework. An important question this poses is: which solutions can be obtained in this way? Or, in other words, can we reconstruct the ambient space field? Since the methods required to address this question are of a different nature to the spirit of this already technical paper, we will address them separately in upcoming work.

\section{A boundary calculus}\label{sec:kg_boundary_calculus}
We return now to the computation of the commutators we postponed in Section~\ref{sec:m_density_space}. These results are of a different nature, and of independent interest from the previous discussion. Indeed, the result we establish in this section constitutes a \emph{boundary calculus} for the operator $H^{AB}D_AD_B$ in our setup. Note that a boundary calculus was established for the \emph{same} operator on Klein-Einstein manifolds in~\cite{smth}, and the setup there is very similar, but we emphasise that we consider a \emph{different} curved orbit decomposition guided by the additional parallel tractor $I^A$.

Boundary calculi have been constructed in the context of Poincaré-Einstein manifolds; we mention in particular the important works of~\cite{Gover:2014aa, gover_poincare-einstein_2015}. The ideas are very much inspired by holographic principles, relating objects on the boundary to objects in the bulk, and in particular the Fefferman-Graham construction. The entire setup of this article, although worded slightly differently, is in exactly this spirit. 

Let us denote by $y$ the operator $H^{AB}D_AD_B$ and $x$ the multiplication operator by the boundary defining density $\sigma$:
\[x: f \mapsto \sigma f,\] both acting on densities of arbitrary weight $w$. The starting point here is to consider their commutator:
\begin{prop}Let $\nu$ be an arbitrary complex number then:
\[ [y,x^\nu]= 2\nu x^{\nu-1}I^AD_A=2\nu x^{\nu-1}D_n. \]
\end{prop}
\begin{proof}
\[\begin{aligned} H^{AB}D_AD_B(x^\nu f)&=H^{AB}\left(fH^{AB}D_AD_Bx^\nu + 2D_{(A}x^\nu D_{B)} f+x^{\nu}D_AD_B f \right)\\&= f\nu (\nu-1)x^{\nu -2}\Phi_{AB}I^AI^B+ \nu x^{\nu-1}D_AI^A + 2\nu x^{\nu-1}I^AD_A f + x^\nu H^{AB}D_AD_B f.\end{aligned} \]
\end{proof}
In particular:
\[ [y,x]=2D_n.\]
Since $\sigma$ is an adapted scale, $D_n\sigma =0$ and therefore $[x,D_n]=0$. We would like to study the operator $[y,D_n]$, we commence with a preliminary computation:
\begin{lemm}\label{prop_commutation_lambda_and_y}
On a section $f$ of $\mathcal{E}(w)$:
\[\begin{aligned} [y,\frac{1}{\sqrt{|\lambda|}}]f &=
&=-\frac{\sgn(\lambda)}{|\lambda|^{\frac{3}{2}}}(n-1+ 2w)f.\end{aligned}\]
\end{lemm}
\begin{proof}
\[\begin{aligned} [y,\frac{1}{\sqrt{|\lambda|}}]f&=H^{AB}D_AD_B(|\lambda|^{-\frac{1}{2}} f) - |\lambda|^{-\frac{1}{2}}H^{AB}D_AD_B f\\
  &=H^{AB}\left(fD_AD_B|\lambda|^{-\frac{1}{2}} + 2D_{(A} |\lambda|^{-\frac{1}{2}}D_{B)}f  \right)\end{aligned}\]
Now:
\[ D_B |\lambda|^{-\frac{1}{2}} = -\frac{1}{2} \frac{\sgn(\lambda)}{|\lambda|^{\frac{3}{2}}} D_B\lambda= -\frac{\sgn(\lambda)}{|\lambda|^{\frac{3}{2}}} \Phi_{BC}X^C. \]
Hence:
\[ D_AD_B|\lambda|^{-\frac{1}{2}}=\frac{3}{|\lambda|^{\frac{5}{2}}}\Phi_{AD}X^D\Phi_{BC}X^C -\frac{\sgn(\lambda)}{|\lambda|^{\frac{3}{2}}}\Phi_{AB},\]
from which it follows that if $f$ is of weight $w$:
\[\begin{aligned} [y,\frac{1}{\sqrt{|\lambda|}}]f &=\frac{\sgn(\lambda)}{|\lambda|^{\frac{3}{2}}}(-(n+2) +3 -2X^AD_A)f\\&=-\frac{\sgn(\lambda)}{|\lambda|^{\frac{3}{2}}}(n-1+ 2w)f.\end{aligned}\]
\end{proof}
In particular, this commutator vanishes if and only if $w=-\frac{n-1}{2}$; this computation shows that one can only hope for the subspaces $\mathcal{E}_m^\pm(w)$ for this specific weight.
Now:

\begin{prop}\label{prop_commutation_of_y_Dn} On densities of arbitrary weight:
\[[y,D_n]=0\]
\end{prop}
\begin{proof}
Note that $D_n f = I^CD_Cf$.
Now, since $I^C$ is parallel we have:
\[H^{AB}D_AD_B I^CD_C= I^CH^{AB}D_AD_BD_C.\]
Now:
\[ \begin{aligned} I^CH^{AB}D_AD_BD_C &= I^CH^{AB}D_A[D_B,D_C] +I^CH^{AB}D_AD_CD_B\\&=I^CH^{AB}D_A[D_B,D_C] + I^CH^{AB}[D_A,D_C]D_B + I^CD_CH^{AB}D_AD_B.\end{aligned}\]
Acting on densities:
\[[D_B,D_C]=0,\]
hence the first term vanishes, and the second term:
\[[D_A,D_C]D_B= (-\Omega_{ac\phantom{E}B}^{\phantom{ac}E}D_E)Z^a_AZ^c_C\equiv -\Omega_{AC\phantom{E}B}^{\phantom{AC}E}D_E,\]
where $\Omega_{ac\phantom{E}B}^{\phantom{ac}E}$ is the curvature tractor.
Overall,
\[[H^{AB}D_AD_B, I^CD_C]=-I^CH^{AB}\Omega_{AC\phantom{E}B}^{\phantom{AC}E}D_E=0.\]
Where the last equality follows, for example, by evaluating in the scale $\sigma$ and using~\cite[Lemma 4.9]{RSTA20230042}.
\end{proof}

This shows that:

\begin{coro}
The Lie algebra generated by $\{x,y\}$ acting on $\mathcal{E}^\pm_m(-\frac{n-1}{2})$ is isomorphic to the Heisenberg algebra,
i.e:
\[ [x,y]=-2D_n, \qquad [x,D_n]=0, \qquad [y,D_n]=0.\]
\end{coro}

%



\newpage

\section{The model case}\label{sec:model}
In this section, we will carry out the theoretical steps in the previous sections in the specific case of the model. This will be done in two ways, first from the projective geometry viewpoint, relying on numerous canonical identifications, and then from the point of view of the base Minkowski spacetime.
\subsection{From the projective sphere}
The starting point here the standard realisation of the projective sphere $S^{n+1}$ as the ray projectivisation of $\R^{n+2}$. The tractor bundle in this case is trivial\footnote{It can be identified with the quotient of $T_{S^{n+1}}$ under the action of the radial vector field}: $\mathcal{T}_{S^{n+1}}=S^{n+1}\times \R^{n+2}$; we identify tractors with vectors in $\R^{n+2}$. For our purposes, we introduce the basis $(e_{\mathbf{A}})$ of $\R^{n+2}$ such that $e_0=\frac{\partial}{\partial X^0} = I$ is identified with the distinguished parallel tractor and the quantity:
\[ h(X,X)= 2X^0X^{n+2} + \eta(Y,Y) \]
where every vector of $\R^{n+2}\setminus\{0\}$ is decomposed as: \[ X= X^0 I +  Y + X^{n+2}e_{n+2}, \quad Y\in \textrm{vect}~(e_{\mathbf{A}})_{\mathbf{A}\in \{1,\dots,n\}}.\]

We emphasise that the vector $X$ represents the homogeneous coordinates of a point on $S^{n+1}$. Introduce now the familiar quantities: \[ \lambda=h(X,X), \qquad \sigma= h(X,I)=X^{n+2}. \]
To begin, we first show how to view things in the scale $\sigma$. In the model, the projection $Z : \mathcal{T}_{S^{n+1}} \longrightarrow TS^{n+1}(-1)$ follows from the identification of functions on $S^{n+1}$ with functions on $\R^{n+2}$ that are homogeneous of degree $0$; in brief one way to think of it is as an interpretation of the derivative of $d\pi$ of the standard projection map.
Trivialising in the scale $\sigma$ is tantamount to working in affine coordinate chart:
\[ X \mapsto \frac{X}{\sigma}=\frac{X}{X^{n+2}}. \]
Restricting to $\Sigma=X^{n+2}>0$, we introduce the coordinate functions $x_0= \frac{X^0}{X^{n+2}}$ and $y^i=\frac{Y^i}{X^{n+2}}$ on $S^{n+1}$, then the expression of the function $l=\lambda \sigma^{-2}$ follows from
\[ \lambda=h(X,X)=(2x_0 + \eta_{ij}y^iy^j)\sigma^2.\] 
Additionally\footnote{we carry out this computation in the scale $\lambda$ with more details below}, we find that:
\[ n=Z(I)=\frac{1}{\sigma}\frac{\partial}{\partial x_0},\quad Z(\frac{\partial}{\partial Y^i})=\frac{1}{\sigma} \frac{\partial}{\partial y^i},\quad Z(\frac{\partial}{\partial X^{n+2}})=\frac{1}{\sigma}\left(-x_0\frac{\partial}{\partial x_0} - y^i\frac{\partial}{\partial y^i}\right).  \]
Applying this to the inverse metric:
\[ H^{AB}=2I^{(A}\frac{\partial^{B)}}{\partial X^{n+2}} + \eta^{ij}\frac{\partial^{A}}{\partial Y^i}\frac{\partial^B}{\partial Y^j},\]
we find the weighted inverse metric:
\[ \zeta^{-1}= \sigma^{-2}\left(2\frac{\partial}{\partial x_0} \left(-x_0\frac{\partial}{\partial x_0} - y^i\frac{\partial}{\partial y^i}\right) + \eta^{ij}\frac{\partial}{\partial y^i}\frac{\partial}{\partial y^j}\right). \]
This looks curious at first, but the coordinate transformation:
\[ \begin{cases} l= 2x_0 + \eta_{ij}y^iy^j ,\\ \tilde{y}^i =y^i, \end{cases} \Rightarrow \begin{cases}\frac{\partial}{\partial x_0} = 2\frac{\partial}{\partial l}, \\ \frac{\partial}{\partial y^i}=\frac{\partial}{\partial\tilde{y}^i} +2\eta_{ij}y^j\frac{\partial}{\partial l},\end{cases} \]
recovers the canonical form expected in the scale $\sigma$ (see~Eq.~\eqref{eq:sigma_scale_metric_decomposition}):
\[\zeta^{ab}=\sigma^{-2}\left(-4l \frac{\partial}{\partial l}\frac{\partial}{\partial l}+ \eta^{ij}\frac{\partial}{\partial \tilde{y}^i}\frac{\partial}{\partial \tilde{y}^j}\right). \]
This recovers explicitly the fact that the interior geometry of the model has a trivial bundle structure.

We now move instead to a boundary scale. One can realise the open set $\lambda \neq 0$ as a union of submanifolds of $\mathbb{R}^{n+2}$: \[ \Lambda=\{ h(X,X)=1\}\cup \{h(X,X)=-1\}\] through the projection:
\[ \pi : \R^{n+2}\setminus\{\lambda =0\} \rightarrow \Lambda,  \]
defined by $\pi(X)=\frac{X}{\sqrt{|h(X,X)|}}.$
This quotients to a homeomorphism\footnote{the inverse map is given by the canonical injection $\Lambda \hookrightarrow \R^{n+2}$ composed with the canonical projection $\R^{n+2} \setminus\{0\} \rightarrow S^{n+1}$} $\frac{\R^{n+2}\setminus\{\lambda =0\}}{\R^*_+}$ onto  $\Lambda$. Now, since when $f$ is a function defined near a point of $\{\lambda \neq 0\}$ in the projective sphere, we can lift it to $\R^{n+2}\setminus\{h(X,X)= 0\}$ using:
\[ \tilde{f}=f\circ \pi,\]
then, the computation of the action of $Z$ on $I$ reads:
\[n(f)= Z(I)(f)= I(\tilde{f})=df_{\pi(X)} \frac{1}{|h(X,X)|^{\frac{1}{2}}}\left(I - \frac{h(X,I)}{h(X,X)}X \right), \]
leading to the weighted vector field $n$. (In the interior scale it was simply $\frac{2}{\sigma}\frac{\partial}{\partial l}$).
Once more, we recognise its projective weight to be $-1$ from its homogeneity with respect to the homogeneous coordinates $X$.
Now, in the boundary scale, it is appropriate to normalise $n^a$ using the scale $\lambda$: 
\[ \hat{N}=\sqrt{|\lambda|}n = I - \frac{h(X,I)}{h(X,X)}X. \]
Following the proof of Proposition~\ref{prop:final_coordinate} we determine its integral curves: 
\[ \gamma'(t) = \hat{N}(\gamma(t)),\]
writing $\gamma(t) = (X_0, y , X_{n+1})$ as before and assuming initial data on \[\Lambda_{\pm}= \{h(X,X)=\pm 1\}\] we can write:
\[\begin{cases} X'_0(t)= 1 - \sgn(h(X,X))X_{n+1}(t)X_0,\\ y'(t)= -\sgn(h(X,X))X_{n+1}(t)y(t), \\ X_{n+2}(t)= - \sgn(h(X,X))X_{n+1}^2(t) \\ h(X(t),X(t))=\sgn(h(X(0),X(0)))\end{cases}.\]
The last equation is immediately solved:
 \[ \frac{1}{X_{n+1}(t)} = \sgn(\lambda)t +  \frac{1}{X_{n+1}(0)} \Rightarrow X_{n+1}(t)= \frac{X_{n+1}(0)}{1+ \sgn(\lambda)X_{n+1}(0) t},\]
 and plugged into the other equations to yield:
 \[ y(t)= \frac{y(0)}{1+\sgn(\lambda)tX_{n+1}(0)}.\]
Finally:
 \[ X_0(t)= \frac{X_0(0)}{1+\sgn(\lambda)tX_{n+1}(0)} + t\frac{1+\frac{1}{2}\sgn(\lambda)tX_{n+1}(0)}{1+t\sgn(\lambda)X_{n+1}(0)}.\]
From here, we can give an explicit parametrisation of a boundary patch given by Prop.~\ref{prop:final_coordinate}. For definiteness, we state the formula on the set $\{ \lambda <0\}$:
\[ (u,\psi,\vartheta,\varrho) \mapsto \frac{1}{1-\varrho u}\begin{pmatrix} u(1-\frac{1}{2}u\varrho)\\ \cosh\psi \\\sinh\psi \vartheta \\\varrho \end{pmatrix}. \]
One can then check that in the scale $\lambda$, the induced metric is given by:
\[ (1-\varrho u)^2\hat{g} = \varrho^2 \textrm{d}u^2 + 2\textrm{d}\varrho \odot \textrm{d}u +\textrm{d}\psi^2 +\sinh^2\psi\textrm{d}\sigma^2,\]
leading to the expression:
\[\bar{g}= \varrho^2 \textrm{d}u^2 + 2\textrm{d}\varrho \odot \textrm{d}u +\underbrace{\textrm{d}\psi^2 +\sinh^2\psi\textrm{d}\sigma^2}_{\bar{h}},\]
which is compatible with Prop.~\ref{prop:final_coordinate}. Furthermore, $s=\sigma|\lambda|^{-\frac{1}{2}}=\frac{\varrho}{1-\varrho u}$, as expected.
One can also check the expression of $\varrho$ in the scale $\sigma$ (and the corresponding coordinates: $(l,y)$):
\[ \begin{cases} \frac{\cosh \psi}{1-\varrho u} = \frac{Y^0}{|\lambda|^{\frac{1}{2}}}=\frac{Y^0}{\sigma}\frac{\varrho}{1-\varrho u} \\ \frac{\sinh\psi}{1-\varrho u}\vartheta^i=\frac{Y^i}{|\lambda|^{\frac{1}{2}}}=\frac{Y^i}{\sigma}\frac{\varrho}{1-\varrho u} \end{cases} \Rightarrow \varrho^2=-\eta_{ij}y^iy^j. \]

An interesting exercise is to compute the expression of the operator $H^{AB}D_AD_B$ in a boundary scale as given by~\eqref{eq:KGinBoundaryScale}. For this, we first reason with the Levi-Civita connection of $\hat{g}$, which is in our working projective class. 
One can check that: 
\[ \textrm{Ric}(\hat{g})_{ab}= -n \hat{g}_{ab}, \]
and therefore the projective Schouten tensor is given by: \[ \hat{P}_{ab}= -\hat{g}_{ab}.\]
Defining the projective scale:
\[\tau= \frac{1}{(1-\varrho u)^2}\lambda,\] we can switch to the connection $\bar{\nabla}$ it preserves in the projective class and find that:
\[ \bar{P}_{ab}=-\bar{h}_{ab}. \]
We can then plug this into Eq.~\eqref{eq:KGinBoundaryScale} giving the following expression for the Klein-Gordon operator in the boundary scale:
\[ -\partial_\varrho(\varrho^2\partial_\varrho) + 2\partial_u\partial_\varrho +\Delta_{\bar{h}} + \frac{(n-1)(n-3)}{4}.  \]

\subsection{From the perspective of Minkowski spacetime}
Rather than view the Minkowski spacetime as a derived structure from a projective geometry, one might instead wish to build the extended spacetime from it.
Start from Minkowski spacetime of dimension $n=d+1$ written in Beig-Schmidt form in the interior of the lightcone.
\begin{equation}
 \eta = -\frac{1}{\varrho^4} d\varrho^2 + \frac{1}{\varrho^2}\bar{h},
\end{equation}
where $\bar{h}$ is the hyperbolic metric of dimension $n-1$; this form is in particular obtained using the coordinate $\varrho=(t^2-|x|^2)^{-\frac{1}{2}}$. Introduce now an extended spacetime, of dimension $n+1$, and a normalised vector field:
\begin{align}\label{eq:extended_minkowski}
\eta' &= - \frac{1}{4\ell}d\ell^2 + \eta,& \hat{N} = 2 \sqrt{-\ell} \partial_{\ell}.
\end{align}
where $\ell<0$. We now introduce the coordinate $u$ such that
\begin{equation}
	\ell = - \frac{(1-\varrho u)^2}{\varrho^2} \quad \Leftrightarrow \quad -u = -\varrho^{-1} + \sqrt{|\ell|}.
\end{equation}
As a result, $\hat{N} = \partial_{u}$ and

\begin{align}
	\eta' &=\frac{1}{\varrho^2} \left( \varrho^2 du^2 + 2 dud\varrho  + \bar{h} \right) ,& \hat{N} =  \partial_{u}.
\end{align}
Consider the conformal rescaling of the metric, recovers the normal form expected in a boundary scale $|\tau|$.
\begin{align}
\overline{\eta}' &=  \varrho^2 du^2 + 2 dud\varrho  + \bar{h},
\end{align}
Here, the boundary of the compactified ambient space, can therefore be seen here to be obtained by a specific \emph{conformal} compactification of the extended metric. 

As in the previous case, we can compute the operator $H^{AB}D_AD_B$ in the underlying boundary scale and coordinates. First we note that if $\phi$ is a projective density of weight $-\frac{n-1}{2}$, then this operator acts on $\varphi = \phi\sigma^\frac{n-1}{2}$, the trivialisation of the function $\varphi$ in the \enquote{physical} scale according to: $\Box^{\eta'}\varphi$, i.e.
\[\sigma^{\frac{n+3}{2}}H^{AB}D_AD_B \phi = \Box^{\eta'}\varphi=-\varrho^4\partial^2_\varrho\varphi +(n-3)\varrho^3\partial_\varrho \varphi +2\varrho^2\partial_u\partial_\varrho\varphi -(n-1)\varrho\partial_u +\varrho^2\Delta_{\bar{h}}\varphi.\]
\begin{rema}
The reader may be concerned that this does not look anything like the Klein-Gordon equation, whereas we expect (by Theorem~\ref{thm:proj_KG}) the operator to impose the Klein-Gordon equation it on the base manifold (assuming the appropriate condition on its behaviour in the direction $\hat{N}$)! The reason is that the boundary coordinates are \emph{not} compatible with the principal bundle structure on the interior, whose global sections correspond to the hypersurface $l=\textrm{cst.}$ However, one can check that if one pulls back $\varphi$ along such a section, we will recognise the Klein-Gordon equation: this also follows from the expression~\eqref{eq:extended_minkowski} of the metric $\eta'$ in the original coordinate system. Indeed:
\[\Box^{\eta'} \varphi= 2\sqrt{-l}(\partial_l(2\sqrt{-l}\partial_l)\varphi + \Box^{\eta} \varphi.\]
\end{rema}

Switching the expression of this operator to the boundary scale $\tau$, amounts to the change of variable 
\[\psi = \varrho^{\frac{n-1}{2}}\varphi,\]
and then:
\begin{equation}\label{eq:KG_boundary_scale_model} \begin{aligned} |\tau|^{\frac{n+3}{2}}H^{AB}D_AD_B\phi &=\varrho^{-\frac{n+3}{2}}\left(\Box^{\eta'}(\psi\varrho^{\frac{n-1}{2}}) \right)\\&=   -\partial_\varrho(\varrho^2\partial_\varrho\psi) + 2\partial_u\partial_\varrho\psi +\Delta_{\bar{h}}\psi + \frac{(n-1)(n-3)}{4}\psi.\end{aligned}\end{equation} 
Recovering the computation of the previous section.

\subsection{A few words on the equation of the boundary field}
We conclude this section with some comments on the content of Section~\ref{sec:KG} in the model case.
The ambient space field $\varphi$ has a canonical decomposition into positive and negative frequency parts, as we have already remarked:
\[ \phi =\varphi\sigma^{-\frac{n-1}{2}}= \varphi_+\sigma^{-\frac{n-1}{2}} + \varphi_-\sigma^{-\frac{n-1}{2}}.  \]
For readability, we will concentrate on $\varphi_+$, but a similar discussion holds for $\varphi_-$. System~\eqref{eq:proj_KGpm} reads here:
\[ \begin{cases} 2\sqrt{|l|}\partial_l\varphi_+ =im\varphi_+, \\  \Box^\eta \varphi_+ - m^2\varphi_+=0.  \end{cases} \]
This is immediately equivalent to the fact that:
\[\varphi_+=\bar{\varphi}_+ e^{-im\sqrt{-l}}, \]
where now: $\bar{\varphi}_+$ is a solution the Klein-Gordon equation on Minkowksi spacetime; i.e. it does not depend on $l$.
On the other hand:
 \[ \varphi_+=\varrho^{\frac{n-1}{2}}\psi_+\]  and $\psi$ is in the kernel of the operator in Eq.~\eqref{eq:KG_boundary_scale_model}.
 In fact: 

\[ \psi_+= \varrho^{-\frac{n-1}{2}}\bar{\varphi}_+e^{-im\sqrt{-l}}= \varrho^{-\frac{n-1}{2}}\bar{\varphi}_+e^{-\frac{im}{\varrho}}e^{imu}\equiv \bar{\psi}e^{imu},\]
where $\partial_u\bar{\psi}=0.$
%

\printbibliography[title=Bibliography]

 \end{document}